\documentclass[15pt, a4paper]{article}
\usepackage{}
\usepackage{amsfonts}
\usepackage{mathbbold}
\usepackage{bbm}
\usepackage{mathrsfs}
\usepackage{latexsym}
\usepackage{graphicx}
\usepackage{amssymb}

\newtheorem{theorem}{Theorem}[section]
\newtheorem{lemma}[theorem]{Lemma}
\newtheorem{corollary}[theorem]{Corollary}

\newtheorem{definition}[theorem]{Definition}

\title{{\Large \bf  On the spectral radii of the unicyclic hypergraphs with fixed matching number\thanks{Supported by NSFC
(No. 11771376, 11571252), Foundation of Lingnan Normal
University(ZL1923), ``333" Project of  Jiangsu (2016), NSFCU of Jiangsu (16KJB110011).}}}
\author{Guanglong Yu$^{a,b}$\thanks{Corresponding author, E-mail addresses:
yglong01@163.com (G. Yu), rockzhang76@tzc.edu.cn (H. Zhang).}
 ~ Chao Yan$^c$ ~ Yarong Wu$^d$ ~  Hailiang Zhang$^e$$^{\dag}$ ~
\\ ~ \\
{\footnotesize $^a$Department of Mathematics, Lingnan Normal
University,  Zhanjiang, 524048, Guangdong, China}\\ {\footnotesize $^b$Department of Mathematics, Yancheng Teachers
University, Yancheng, 224002, Jiangsu, China}\\
{\footnotesize $^c$Department of Mathematics, Pujiang Institute, Nanjing Tech University, Nanjing, 211101, Jiangsu, China}\\
{\footnotesize $^d$ SMU college of art and science, Shanghai maritime
University, Shanghai, 200135, China}\\
{\footnotesize $^e$Department of Mathematics, Taizhou University, Linhai, Zhejiang, 317000, China}}
\date{}

\begin{document}
\maketitle

\begin{abstract}
We determine the unique hypergraphs with maximum spectral radius among all connected $k$-uniform ($k\geq 3$) unicyclic hypergraphs with matching number at least $z$, and among all connected $k$-uniform ($k\geq 3$) unicyclic hypergraphs with a given matching number, respectively.


\end{abstract}
\baselineskip 18.6pt

\section{Introduction}

\ \ \ \ In the 19th
century, Gauss et al. had introduced the concept of tensors in the study of differential geometry. In the very beginning of the 20th century, Ricci et al. further developed tensor analysis as a mathematical discipline. It was Einstein
who applied tensor analysis in his study of general relativity in 1916. This made tensor analysis become an important tool in theoretical physics, continuum mechanics and many other areas of science
and engineering. Thus a comprehensive
study of issues relevant to tensors has been undertaken \cite{DGAO, KCHX, KTLH, LLS, LQIE, YRTT}.

Denote by $\mathbb{C}$ the set of
all complex numbers, $\mathbb{R}$ the set of
all real numbers, $\mathbb{R}_{+}$ the set of
all nonegative real numbers and $\mathbb{R}^{n}_{++}$ the set of all positive real numbers. Recall that a $k$th-order $n$-dimensional real tensor $\mathcal {T}$ consists of $n^{k}$ entries in real numbers:
$\mathcal {T} = (\mathcal {T}_{i_{1}\cdots i_{k}}), \mathcal {T}_{i_{1}\cdots i_{k}}\in R$, $1 \leq i_{1}, \ldots, i_{k} \leq n$. $\mathcal {T}$ is called symmetric if the value of $\mathcal {T}_{i_{1}\cdots i_{k}}$ is invariant under any permutation of
its indices $i_{1}, \ldots, i_{k}$. $\mathcal {T}X^{k-1}$ is a vector in $\mathbb{C}^{n}$ with its ith component as
$(\mathcal {T}X^{k-1})_{i} =\sum^{n}_{i_{2}, \ldots, i_{k}=1}
\mathcal {T}_{i, i_{2}, \ldots, i_{k}} x_{i_{2}}\cdots x_{i_{k}}$, where $X=(x_1, x_2, \ldots, x_n)^T\in \mathbb{C}^{n}$; for $X=(x_1, x_2, \ldots, x_n)^T\in \mathbb{R}^{n}$, $X^{T}\mathcal {T}X^{k-1}=\sum^{n}_{i_{1}, i_{2}, \ldots, i_{k}=1}
\mathcal {T}_{i_{1}, i_{2}, \ldots, i_{k}} x_{i_{1}}x_{i_{2}}\cdots x_{i_{k}}$.

\begin{definition}{\bf \cite{KCKPZ}, \cite{LQIE}}\label{de01,01} 
Let $\mathcal {T}$ be a kth-order n-dimensional tensor. A pair $(\lambda, X)$ $(X\in\mathbb{C}^{n}\setminus \{0\})$ is called an eigenvalue and an eigenvector of $\mathcal {T}$ if they satisfy
$\mathcal {T}X^{k-1}=\lambda X^{[k-1]} $, where $X^{[k-1]}=(x^{k-1}_{1}, x^{k-1}_{2}, \ldots, x^{k-1}_{n})^T$.
\end{definition}

The $spectrum$ of a real symmetric tensor $\mathcal {T}$ is defined as the multiset of its eigenvalues, and its $spectral$ $radius$, denoted by $\rho(\mathcal {T})$, is the maximum modulus among its all eigenvalues.

\begin{lemma}{\bf \cite{LQI}}\label{le01.02} 
Let $\mathcal {T}$ be a kth-order $n$-dimensional nonnegative symmetric tensor. Then we have
$\rho(\mathcal {T}) = \max\{X^{T}\mathcal {T}X^{k-1}\mid
\sum^{n}_{i=1}x^{k}_{i}= 1, X\in R^{n}_{+}\}$. Furthermore, $x\in R^{n}_{+}$ with $\sum^{n}_{i=1}x^{k}_{i}=1$ is an optimal solution of above optimization problem if and only if it is an eigenvector of $\mathcal {T}$ corresponding to the eigenvalue $\rho(\mathcal {T})$.

\end{lemma}

In recent years, since the work of Qi \cite{LQIE} and Lim \cite{LLS}, the study of the spectra of tensors and hypergraphs with their various applications has attracted extensive attention and interest. In 2008, Lim \cite{LLE} proposed the study of the spectra of hypergraphs via the spectra of tensors. In 2012, Cooper and Dutle \cite{CJDA} defined the spectrum of a uniform hypergraph as the spectrum of the adjacency tensor of that hypergraph, and obtained hypergraph generalizations of many basic results of spectral graph theory. The (adjacency) spectrum of uniform hypergraphs were further studied in \cite{MFYZ, LKNY, LSQ, JZLS, XYJS}.

Recall that a $hypergraph$ $G=(V, E)$ consists of a vertex set $V=V(G)$ and an edge set $E=E(G)$, where $V=V(G)$ is nonempty and each edge $e\in E(G)$ is a subset of $V(G)$ containing at least two elements. The cardinality $n=\|V(G)\|$ is called the order; $m=\|E(G)\|$ is called the edge number of hypergraph $G$. Denote by $t$-set a set with size (cardinality) $t$. We say that a hypergraph $G$ is $uniform$ if its every edge has the same size, and call it $k$-$uniform$ if its every edge has size $k$ (i.e. every edge is a $k$-subset). It is well known that a $2$-graph is the general graph. The $adjacency$ $tensor$ $\mathcal {A} =\mathcal {A}(G)$ of a $k$-uniform graph $G$ refers to a multi-dimensional array with entries $\mathcal {A}_{i_{1}\cdots i_{k}}$ such that
$$\mathcal {A}_{i_{1}\cdots i_{k}}= \left \{\begin{array}{ll}
               \frac{1}{(k-1)! }, \ & \ if\ \{i_{1}, \ldots, i_{k}\}\ is\ an\ edge\ of\ G\\
               0, \ & \ others.
             \end{array}\right.$$
where each $i_{j}$ runs from 1 to $n$ for $j\in[k]$ $([k]=\{1, 2, \ldots, k\})$. It can be seen that $\mathcal {A}$ is symmetric. For convenience, the spectrum of the adjacency tensor of hypergraph $G$ is called the spectrum of $G$, the spectral radius $\rho(\mathcal {A})$ is called the spectral radius of $G$. We employ $\rho(G)$ to denote the spectral radius of $G$ without discrimination.
From Lemma \ref{le01.02}, we know that for a $k$-uniform hypergraph $G$, $\rho(G)$ is the optimization of the system $f(X)=X^{T}\mathcal {A}X^{k-1}$ based on this hypergraph under the condition $\sum^{n}_{i=1}x^{k}_{i}= 1, X\in R^{n}_{+}$. In spectral theory of hypergraphs, the spectral radius is an index that attracts much attention due to the fine properties \cite{KCKPZ, YFTP, SFSH, LSQ, HLBZ, LLSM, COQY, YQY}.

We assume that the hypergraphs throughout
the paper are simple, i.e. $e_{i} \neq e_{j}$ if $i \neq j$. For a hypergraph $G$, we define $G-e$ ($G+e$)
to be the graph obtained from $G$ by deleting the edge $e\in
 E(G)$ (by adding an new edge $e$ if $e\notin
 E(G)$); for a edge subset $B\subseteq E(G)$, we define $G-B$
to be the graph obtained from $G$ by deleting each edge $e\in
 B$. For two $k$-uniform hypergraphs $G_{1}=(V_{1}, E_{1})$ and $G_{2}=(V_{2}, E_{2})$, we say the two graphs are $isomorphic$ if there is a bijection $f$ from $V_{1}$ to $V_{2}$, and there is a bijection $g$ from $E_{1}$ to $E_{2}$ that maps each edge $\{v_{1}$, $v_{2}$, $\ldots$, $v_{k}\}$ to $\{f(v_{1})$, $f(v_{2})$, $\ldots$, $f(v_{k})\}$.

In a hypergraph, two vertices are said to be $adjacent$ if
both of them is contained in an edge.  The $neighbor$ $set$ of vertex $v$ in hypergraph $G$, denoted by $N_{G}(v)$, is the set of vertices adjacent to $v$ in $G$. Two edges are said to be $adjacent$
if their intersection is not empty. An edge $e$ is said to be $incident$ with a vertex $v$ if
$v\in e$. The $degree$ of a vertex $v$ in $G$, denoted by $deg_{G}(v)$ (or $deg(v)$ for short), is the number of the edges incident with $v$. For a hypergraph $G$, among all its vertices, we denote by $\Delta(G)$ (or $\Delta$ for short) the $maximal$ $degree$, and denote by $\delta(G)$ (or $\delta$ for short) the $minimal$ $degree$ respectively. A vertex of degree $1$ is called a $pendant$ $vertex$. A $pendant$ $edge$ is an edge with exactly one vertex of degree more than one and other vertices in this edge being all pendant vertices. A pendant vertex in a pendant edge is called a $PP$-$vertex$ here, and a vertex which is not a $PP$-vertex is called a $NPP$-$vertex$.  In a hypergraph $G$, a $path$ of length $q$ ($q$-$path$) is defined to be an alternating sequence
of vertices and edges $v_{1}e_{1}v_{2}e_{2}\cdots v_{q}e_{q}v_{q+1}$ such that
(1) $v_{1}$, $v_{2}$, $\ldots$, $v_{q+1}$ are all distinct vertices;
(2) $e_{1}$, $e_{2}$, $\ldots$, $e_{q}$ are all distinct edges;
(3) $v_{i}$, $v_{i+1}\in e_{i}$ for $i = 1$, $2$, $\ldots$, $q$.
A $cycle$ of length $q$ ($q$-$cycle$) $v_{1}e_{1}v_{2}e_{2}\cdots v_{q}e_{q}v_{1}$ is obtained from a path $v_{1}e_{1}v_{2}e_{2}\cdots v_{q}$ by adding a new edge $e_{q}$ between $v_{1}$ and $v_{q}$ where $e_{q}\neq e_{i}$ for $1\leq i\leq q-1$. We denote by $L(C)$ the length of a cycle $C$. A hypergraph $G$
is connected if there exists a path starting at $v$ and terminating at $u$ for all $v, u \in V$,
and $G$ is called $acyclic$ if it contains no cycle. In a connected hypergraph $G$, the $distance$ of vertex $u, v$, denoted by $d(u, v)$ or $d_{G}(u, v)$, is the length of the shortest path from $u$ to $v$.

A hypergraph $G$ is called a $linear$ $hypergraph$, if each pair of the
edges of $G$ have at most one common vertex. A $supertree$ is a hypergraph which is both connected and acyclic.
Clearly, an acyclic hypergraph is linear. A connected $k$-uniform hypergraph with $n$ vertices and $m$ edges is $r$-cyclic if $n-1 =(k-1)m-r$. For $r=1$ or $2$, it is called a $unicyclic$ $hypergraph$ or a $bicyclic$ $hypergraph$; for $r=0$, it is a supertree.
In \cite{COQY}, C. Ouyang et al. proved that a simple connected $k$-graph $G$ is unicyclic (1-cyclic) if and only if it has only one cycle. From this, for a unicyclic hypergraph $G$ with unique cycle $C$, it follows that (1) if $L(C)=2$, then the two edges in $C$ have exactly two common vertices, and $\|e\cap f\|\leq 1$ for any two edges $e$ and $f$ not in $C$ simultaneously; (2) if $L(C)\geq 3$, then any two edges in $G$ have at most one common vertices; (3) every connected component of $G-E(C)$ is a supertree.

From \cite{SFSH}, for a connected uniform hypergraph $G$ of order $n$, we know that there is unique one positive eigenvector $X=(x(v_{1})$, $x(v_{2})$, $\ldots$, $x(v_{n}))^T \in R^{n}_{++}$ corresponding to $\rho(G)$ where $\sum^{n}_{i=1}x^{k}(v_{i})= 1$,  each vertex $v_i$ is mapped to $x(v_i)$. We call such an eigenvector $X$
the $principal$ $eigenvector$ of $G$. In the principal eigenvector, a vertex $u$ is said to be a $M_{a}$-$vertex$ if $x(u)=\max\{x(v)\mid v\in V(G)\}$.

A $matching$ of hypergraph $G$ is a subset of
independent (pairwise nonadjacent) edges of $E(G)$. The \emph{matching number} of $G$, denoted
by $\alpha(G)$, is the maximum of the cardinalities of all matchings. A $maximal$ $matching$ of $G$ is a matching of $G$ with cardinality $\alpha(G)$. The topic about matching and matching number of a graph or a hypergraph always attract the researchers. In \cite{YHJL}, the authors determined the unique general tree with maximum spectral radius among all connected general trees with a fixed matching number. In \cite{AYFT}, the authors determined the unique unicyclic  general graphs with maximum spectral radius among all the unicyclic general graphs with with fixed matching number. Recently, in \cite{HGBZ}, the authors determined the unique supertrees with maximum spectral radius among all connected $k$-uniform supertrees with a fixed matching number. Note that from the acyclic graph to cyclic graph, the idea or techniques of researching problems always need some quite different transits. Now, an natural problem arising is that how about the maximum spectral radius among all connected $k$-uniform ($k\geq 3$) cyclic hypergraphs with a fixed matching number. Motivated by this, we consider the maximum spectral radius among all connected $k$-uniform ($k\geq 3$) unicyclic hypergraph with a fixed matching number.

Let $\mathscr{U}(n,k;f;r,s;t,w)$ ($f\leq 1$, $r\leq k-2$, $s\leq k-2$, $w\leq k-2$) be a $k$-uniform unicyclic hypergraph of order $n$ obtained from a 2-cycle $C=v_{1}e_{1}v_{2}e_{2}v_{1}$ by: (1) attaching $f$ pendant edge to vertex $v_{2}$; (2) attaching $r$ pendant edges to $r$ vertices of $e_{1}\setminus\{v_{1}, v_{2}\}$ respectively (i.e. each of these $r$ vertices is attached exactly one pendant edge); (3) attaching $s$ pendant edges to $s$ vertices of $e_{2}\setminus\{v_{1}, v_{2}\}$ respectively; (4) attaching $t$ nonpendant edge to $v_{1}$ which is neither $e_{1}$ nor $e_{2}$, where for every edge of these $t$ edges, each vertex of this edge other than $v_{1}$ is attached exactly one pendant edge; (5) attaching one nonpendant edges to $v_{1}$ which is neither $e_{1}$ nor $e_{2}$, where there is a $w$-subset (a subset with cardinality $w$) of this edge containing no $v_{1}$ that each vertex of this subset is attached exactly one pendant edge; (6) attaching $m-f-r-s-t-t(k-1)-w-\lceil\frac{w}{k-1}\rceil-2$ pendant edges to $v_{1}$, where $m=\frac{n}{k-1}$ (For example, see figures shown in Fig. 1.1).

Denote by $\mathcal {H}=\{G\mid G$ is a $k$-uniform unicyclic hypergraph of order $n$ and with $\alpha(G)\geq z$, where $k\geq 3\}$, $\mathbb{H}=\{G\mid G\in \mathcal {H}$ and $\alpha(G)= z\}$, $\rho_{max}=\max\{\rho(G)\mid G\in \mathcal {H}\}$, $\rho^{\ast}_{max}=\max\{\rho(G)\mid G\in \mathbb{H}\}$. In this paper, we determine the $\rho_{max}$ and $\rho^{\ast}_{max}$, getting the following results.

\setlength{\unitlength}{0.5pt}
\begin{center}
\begin{picture}(873,220)
\qbezier(0,141)(0,156)(4,167)\qbezier(4,167)(9,178)(17,178)\qbezier(17,178)(24,178)(29,167)\qbezier(29,167)(34,156)(34,141)\qbezier(34,141)(34,125)(29,114)\qbezier(29,114)(24,103)(17,103)
\qbezier(17,104)(9,104)(4,114)\qbezier(4,114)(0,125)(0,141)
\put(16,177){\circle*{4}}
\put(17,105){\circle*{4}}
\qbezier(16,177)(0,125)(17,105)
\qbezier(16,177)(32,126)(17,105)
\qbezier(23,169)(23,167)(29,166)\qbezier(29,166)(36,165)(47,165)\qbezier(47,165)(57,165)(64,166)\qbezier(64,166)(71,167)(71,169)\qbezier(71,169)(71,170)(64,171)\qbezier(64,171)(57,172)(47,172)
\qbezier(47,172)(36,172)(29,171)\qbezier(29,171)(22,170)(23,169)
\qbezier(26,153)(26,154)(33,155)\qbezier(33,155)(40,156)(50,156)\qbezier(50,156)(59,156)(66,155)\qbezier(66,155)(74,154)(74,153)\qbezier(74,153)(74,151)(66,150)\qbezier(66,150)(59,149)(50,149)
\qbezier(50,150)(40,150)(33,150)\qbezier(33,150)(26,151)(26,153)
\qbezier(28,118)(28,119)(34,120)\qbezier(34,120)(41,121)(51,121)\qbezier(51,121)(60,121)(67,120)\qbezier(67,120)(74,119)(74,118)\qbezier(74,118)(74,116)(67,115)\qbezier(67,115)(60,114)(51,114)
\qbezier(51,115)(41,115)(34,115)\qbezier(34,115)(27,116)(28,118)
\put(60,142){\circle*{4}}
\put(60,135){\circle*{4}}
\put(60,127){\circle*{4}}
\qbezier(263,153)(263,168)(258,179)\qbezier(258,179)(252,190)(245,190)\qbezier(245,190)(237,190)(231,179)\qbezier(231,179)(227,168)(226,153)\qbezier(226,153)(227,137)(231,126)\qbezier(231,126)(237,116)(245,116)
\qbezier(245,116)(252,116)(258,126)\qbezier(258,126)(263,137)(263,153)
\put(245,190){\circle*{4}}
\put(245,117){\circle*{4}}
\qbezier(245,190)(228,144)(245,117)
\qbezier(245,190)(260,142)(245,117)
\qbezier(244,190)(244,191)(251,193)\qbezier(251,193)(259,194)(270,194)\qbezier(270,194)(280,194)(288,193)\qbezier(288,193)(296,191)(296,190)\qbezier(296,190)(296,188)(288,186)\qbezier(288,186)(280,185)(270,185)
\qbezier(270,186)(259,186)(251,186)\qbezier(251,186)(244,188)(244,190)
\put(143,193){\circle*{4}}
\put(144,118){\circle*{4}}
\put(726,180){\circle*{4}}
\qbezier(162,154)(162,169)(157,181)\qbezier(157,181)(151,192)(144,192)\qbezier(144,192)(136,192)(130,181)\qbezier(130,181)(126,169)(125,154)\qbezier(125,154)(126,138)(130,126)\qbezier(130,126)(136,115)(144,115)
\qbezier(144,116)(151,116)(157,126)\qbezier(157,126)(162,138)(162,154)
\qbezier(143,193)(128,140)(144,118)
\qbezier(143,193)(157,137)(144,118)
\put(137,100){\circle*{4}}
\qbezier(144,118)(167,115)(175,101)
\qbezier(144,118)(162,96)(175,101)
\qbezier(144,118)(118,114)(115,100)
\qbezier(144,118)(125,97)(115,100)
\put(144,100){\circle*{4}}
\put(151,100){\circle*{4}}
\qbezier(245,117)(275,114)(278,102)
\qbezier(245,117)(262,100)(278,102)
\qbezier(245,117)(217,114)(216,101)
\qbezier(216,101)(226,98)(245,117)
\put(246,101){\circle*{4}}
\put(238,101){\circle*{4}}
\put(254,101){\circle*{4}}
\put(363,192){\circle*{4}}
\put(363,119){\circle*{4}}
\put(334,103){\circle*{4}}
\put(364,103){\circle*{4}}
\put(356,103){\circle*{4}}
\put(372,103){\circle*{4}}
\qbezier(381,155)(381,170)(376,181)\qbezier(376,181)(370,192)(363,192)\qbezier(363,192)(355,192)(349,181)\qbezier(349,181)(345,170)(344,155)\qbezier(344,155)(345,139)(349,128)\qbezier(349,128)(355,118)(363,118)
\qbezier(363,118)(370,118)(376,128)\qbezier(376,128)(381,139)(381,155)
\qbezier(363,192)(346,146)(363,119)
\qbezier(363,192)(378,144)(363,119)
\qbezier(362,192)(362,193)(369,195)\qbezier(369,195)(377,196)(388,196)\qbezier(388,196)(398,196)(406,195)\qbezier(406,195)(414,193)(414,192)\qbezier(414,192)(414,190)(406,188)\qbezier(406,188)(398,187)(388,187)
\qbezier(388,188)(377,188)(369,188)\qbezier(369,188)(362,190)(362,192)
\qbezier(363,119)(393,116)(396,104)
\qbezier(363,119)(380,102)(396,104)
\qbezier(363,119)(335,116)(334,103)
\qbezier(334,103)(344,100)(363,119)
\put(407,162){\circle*{4}}
\put(407,155){\circle*{4}}
\put(407,147){\circle*{4}}
\qbezier(376,173)(376,174)(383,175)\qbezier(383,175)(390,176)(400,176)\qbezier(400,176)(409,176)(416,175)\qbezier(416,175)(424,174)(424,173)\qbezier(424,173)(424,171)(416,170)\qbezier(416,170)(409,169)(400,169)
\qbezier(400,170)(390,170)(383,170)\qbezier(383,170)(376,171)(376,173)
\qbezier(375,138)(375,139)(382,140)\qbezier(382,140)(389,141)(400,141)\qbezier(400,141)(410,141)(417,140)\qbezier(417,140)(425,139)(425,138)\qbezier(425,138)(425,136)(417,135)\qbezier(417,135)(410,134)(400,134)
\qbezier(400,135)(389,135)(382,135)\qbezier(382,135)(375,136)(375,138)
\put(519,189){\circle*{4}}
\put(519,116){\circle*{4}}
\put(520,100){\circle*{4}}
\put(512,100){\circle*{4}}
\put(528,100){\circle*{4}}
\put(563,159){\circle*{4}}
\put(563,152){\circle*{4}}
\put(563,144){\circle*{4}}
\qbezier(537,152)(537,167)(532,178)\qbezier(532,178)(526,189)(519,189)\qbezier(519,189)(511,189)(505,178)\qbezier(505,178)(501,167)(500,152)\qbezier(500,152)(501,136)(505,125)\qbezier(505,125)(511,115)(519,115)
\qbezier(519,115)(526,115)(532,125)\qbezier(532,125)(537,136)(537,152)
\qbezier(519,189)(502,143)(519,116)
\qbezier(519,189)(534,141)(519,116)
\qbezier(518,189)(518,190)(525,192)\qbezier(525,192)(533,193)(544,193)\qbezier(544,193)(554,193)(562,192)\qbezier(562,192)(570,190)(570,189)\qbezier(570,189)(570,187)(562,185)\qbezier(562,185)(554,184)(544,184)
\qbezier(544,185)(533,185)(525,185)\qbezier(525,185)(518,187)(518,189)
\qbezier(519,116)(549,113)(552,101)
\qbezier(519,116)(536,99)(552,101)
\qbezier(519,116)(491,113)(490,100)
\qbezier(490,100)(500,97)(519,116)
\qbezier(532,170)(532,171)(539,172)\qbezier(539,172)(546,173)(556,173)\qbezier(556,173)(565,173)(572,172)\qbezier(572,172)(580,171)(580,170)\qbezier(580,170)(580,168)(572,167)\qbezier(572,167)(565,166)(556,166)
\qbezier(556,167)(546,167)(539,167)\qbezier(539,167)(532,168)(532,170)
\qbezier(531,135)(531,136)(538,137)\qbezier(538,137)(545,138)(556,138)\qbezier(556,138)(566,138)(573,137)\qbezier(573,137)(581,136)(581,135)\qbezier(581,135)(581,133)(573,132)\qbezier(573,132)(566,131)(556,131)
\qbezier(556,132)(545,132)(538,132)\qbezier(538,132)(531,133)(531,135)
\put(473,145){\circle*{4}}
\put(473,153){\circle*{4}}
\put(473,160){\circle*{4}}
\qbezier(458,136)(458,137)(465,138)\qbezier(465,138)(472,139)(483,139)\qbezier(483,139)(493,139)(500,138)\qbezier(500,138)(508,137)(508,136)\qbezier(508,136)(508,134)(500,133)\qbezier(500,133)(493,132)(483,132)
\qbezier(483,133)(472,133)(465,133)\qbezier(465,133)(458,134)(458,136)
\qbezier(460,171)(460,172)(467,173)\qbezier(467,173)(474,174)(484,174)\qbezier(484,174)(493,174)(500,173)\qbezier(500,173)(508,172)(508,171)\qbezier(508,171)(508,169)(500,168)\qbezier(500,168)(493,167)(484,167)
\qbezier(484,168)(474,168)(467,168)\qbezier(467,168)(460,169)(460,171)
\put(752,135){\circle*{4}}
\qbezier(627,139)(693,150)(752,135)
\qbezier(627,139)(689,120)(752,135)
\qbezier(627,160)(627,151)(628,145)\qbezier(628,145)(629,139)(631,139)\qbezier(631,139)(632,139)(633,145)\qbezier(633,145)(635,151)(635,160)\qbezier(635,160)(635,168)(633,174)\qbezier(633,174)(632,181)(631,181)
\qbezier(631,181)(629,181)(628,174)\qbezier(628,174)(627,168)(627,160)
\qbezier(650,166)(650,157)(650,151)\qbezier(650,151)(651,145)(653,145)\qbezier(653,145)(654,145)(655,151)\qbezier(655,151)(656,157)(656,166)\qbezier(656,166)(656,174)(655,180)\qbezier(655,180)(654,187)(653,187)
\qbezier(653,187)(651,187)(650,180)\qbezier(650,180)(649,174)(650,166)
\qbezier(694,166)(694,173)(694,179)\qbezier(694,179)(694,184)(696,184)\qbezier(696,184)(697,184)(697,179)\qbezier(697,179)(698,173)(698,166)\qbezier(698,166)(698,158)(697,152)\qbezier(697,152)(697,147)(696,147)
\qbezier(696,148)(694,148)(694,152)\qbezier(694,152)(693,158)(694,166)
\put(687,104){\circle*{4}}
\put(672,104){\circle*{4}}
\put(680,104){\circle*{4}}
\put(682,164){\circle*{4}}
\put(667,164){\circle*{4}}
\put(675,164){\circle*{4}}
\put(743,190){\circle*{4}}
\qbezier(752,135)(734,165)(743,190)
\qbezier(752,135)(758,168)(743,190)
\qbezier(743,190)(722,203)(714,191)
\qbezier(743,190)(720,184)(714,191)
\qbezier(748,154)(728,149)(719,161)
\qbezier(719,161)(724,163)(748,154)
\put(726,167){\circle*{4}}
\put(726,174){\circle*{4}}
\put(785,176){\circle*{4}}
\qbezier(752,135)(783,156)(785,176)
\qbezier(752,135)(768,178)(785,176)
\put(765,110){\circle*{4}}
\put(785,177){\circle*{4}}
\qbezier(821,187)(802,192)(785,177)
\put(785,177){\circle*{4}}
\qbezier(823,187)(812,171)(785,177)
\qbezier(809,151)(796,145)(762,149)
\qbezier(762,149)(795,162)(809,152)
\qbezier(733,203)(768,219)(805,197)
\put(767,221){$t$}
\put(833,137){\circle*{4}}
\qbezier(752,135)(789,138)(833,137)
\qbezier(752,135)(838,118)(833,137)
\put(815,109){\circle*{4}}
\qbezier(779,133)(787,111)(815,109)
\qbezier(779,133)(810,119)(815,109)
\put(867,124){\circle*{4}}
\qbezier(833,137)(861,118)(867,124)
\qbezier(833,137)(867,134)(867,124)
\qbezier(813,103)(866,101)(873,117)
\put(854,94){$w$}
\put(759,108){\circle*{4}}
\put(752,107){\circle*{4}}
\put(824,115){\circle*{4}}
\put(831,117){\circle*{4}}
\put(837,120){\circle*{4}}
\qbezier(628,136)(628,143)(646,149)\qbezier(646,149)(664,155)(690,155)\qbezier(690,155)(715,155)(733,149)\qbezier(733,149)(752,143)(752,136)\qbezier(752,136)(752,128)(733,122)\qbezier(733,122)(715,117)(690,117)
\qbezier(690,117)(664,117)(646,122)\qbezier(646,122)(628,128)(628,136)
\qbezier(658,125)(634,106)(640,94)
\qbezier(658,125)(654,100)(640,94)
\qbezier(709,124)(699,105)(699,94)
\qbezier(709,124)(713,97)(699,94)
\qbezier(752,135)(771,124)(780,110)
\qbezier(752,135)(767,112)(780,110)
\qbezier(752,135)(733,119)(735,101)
\qbezier(752,135)(748,107)(735,101)
\put(807,170){\circle*{4}}
\put(803,166){\circle*{4}}
\put(799,162){\circle*{4}}
\put(630,139){\circle*{4}}
\put(8,91){$v_{1}$}
\put(10,187){$v_{2}$}
\put(138,125){$v_{1}$}
\put(138,198){$v_{2}$}
\put(239,125){$v_{1}$}
\put(231,196){$v_{2}$}
\put(358,125){$v_{1}$}
\put(353,201){$v_{2}$}
\put(514,124){$v_{1}$}
\put(508,198){$v_{2}$}
\put(728,133){$v_{1}$}
\put(607,138){$v_{2}$}
\put(9,48){$G_{1}$}
\put(138,48){$G_{2}$}
\put(237,48){$G_{3}$}
\put(357,48){$G_{4}$}
\put(511,48){$G_{5}$}
\put(729,48){$G_{6}$}
\put(347,3){Fig. 1.1. $G_{1}-G_{6}$}
\put(769,177){\circle*{4}}
\put(754,184){\circle*{4}}
\put(762,181){\circle*{4}}
\end{picture}
\end{center}

\begin{theorem}\label{th01.03} 
Let $G\in \mathcal {H}$ and $\rho(G)=\rho_{max}$. Then $m\geq z+1$, and

$\mathrm{(1)}$ if $m=z+1$, then $G\cong  \mathscr{U}(n,k;0;0,z-1;0,0)$ (see $G_{1}$ in Fig. 1.1);

$\mathrm{(2)}$ if $m\geq z+2$, $1\leq z\leq 2$, then $G\cong  \mathscr{U}(n,k;z-1;0,0;0,0)$ (see $G_{2}$, $G_{3}$ in Fig. 1.1);

$\mathrm{(3)}$ if $m\geq z+2$, $3\leq z\leq k$, then $G\cong \mathscr{U}(n,k;1;z-2,0;0,0)$ (see $G_{4}$ in Fig. 1.1);

$\mathrm{(4)}$ if $m\geq z+2$, $k+1\leq z\leq 2k-2$, then $G\cong \mathscr{U}(n,k;1;k-2,z-k;0,0)$ (see $G_{5}$ in Fig. 1.1);

$\mathrm{(5)}$ if $m\geq z+2$, $2k-1\leq z\leq m-2-\lceil\frac{m-2k}{k-1}\rceil$, then $G\cong \mathscr{U}(n,k;1;k-2,k-2;t, w)$, where $t=\lfloor\frac{z-2k+2}{k-1}\rfloor$, $w=z-2k+2-\lfloor\frac{z-2k+2}{k-1}\rfloor(k-1)$ (see $G_{6}$ in Fig. 1.1).
\end{theorem}

\begin{theorem}\label{th01.04} 
Let $G\in \mathbb{H}$ and $\rho(G)=\rho^{\ast}_{max}$. Then $m\geq z+1$, and

$\mathrm{(1)}$ if $m=z+1$, then $G\cong  \mathscr{U}(n,k;0;0,z-1;0,0)$ (see $G_{1}$ in Fig. 1.1);

$\mathrm{(2)}$ if $m\geq z+2$, $1\leq z\leq 2$, then $G\cong  \mathscr{U}(n,k;z-1;0,0;0,0)$ (see $G_{2}$, $G_{3}$ in Fig. 1.1);

$\mathrm{(3)}$ if $m\geq z+2$, $3\leq z\leq k$, then $G\cong \mathscr{U}(n,k;1;z-2,0;0,0)$ (see $G_{4}$ in Fig. 1.1);

$\mathrm{(4)}$ if $m\geq z+2$, $k+1\leq z\leq 2k-2$, then $G\cong \mathscr{U}(n,k;1;k-2,z-k;0,0)$ (see $G_{5}$ in Fig. 1.1);

$\mathrm{(5)}$ if $m\geq z+2$, $2k-1\leq z\leq m-2-\lceil\frac{m-2k}{k-1}\rceil$, then $G\cong \mathscr{U}(n,k;1;k-2,k-2;t, w)$, where $t=\lfloor\frac{z-2k+2}{k-1}\rfloor$, $w=z-2k+2-\lfloor\frac{z-2k+2}{k-1}\rfloor(k-1)$ (see $G_{6}$ in Fig. 1.1).
\end{theorem}

The layout of this paper is as follows: section 2 introduces some notations and working lemmas; section 3 represents the main results.

\section{Preliminary}

\ \ \ \ \ In this section, we introduce some notations and some working lemmas.

From \cite{JYS}, for two tensors $\mathcal {D}, \mathcal {T}$, we know that if there exists a permutation matrix $P = P_{\sigma}$ (corresponding to a permutation $\sigma\in S_{n}$) such that
$\mathcal {D} = P\mathcal {T}P^{T}$, where $\mathcal {D}_{i_{1},\ldots,i_{k}} = \mathcal {T}_{\sigma(i_{1}),\sigma(i_{2}),\ldots,\sigma(i_{k})}$ then $\mathcal {D}$ and $\mathcal {T}$ are called permutational similar.
As the proof of Proposition 27 and Lemma 28 in \cite{LSQ} (change $\mathcal {Q}^{\ast}$ into $\mathcal {A}$), we can get the following two results.

\begin{lemma}\label{le02,01} 
Let $G$ be a $k$-uniform hypergraph on $n$ vertices, $\mathcal {A}$
be its
adjacency tensor. A permutation $\sigma \in S_{n}$ is an automorphism of
$G$ if and only if $P_{\sigma}\mathcal {A} =\mathcal {A}P_{\sigma}$, where $P_{\sigma}$ is a permutation matrix corresponding to a permutation $\sigma\in S_{n}$.
\end{lemma}

\begin{lemma}\label{le02,02} 
Let $G$ be a connected $k$-uniform hypergraph, $\mathcal {A}$
be its
adjacency tensor. If $X$ is an eigenvector of $\mathcal {A}$ corresponding to
eigenvalue $\lambda$, then for each automorphism $\sigma$ of $G$, we have $P_{\sigma} X$ is also an eigenvector of $\mathcal {A}$ corresponding to the eigenvalue $\lambda$.
Moreover, if $X$ is an eigenvector of $\mathcal {A}$ corresponding to
eigenvalue $\rho(\mathcal {A})$, then we have

$\mathrm{(1)}$ $P_{\sigma} X = X$;

$\mathrm{(2)}$ For any orbit $\Omega$ of $Aut (G)$ and each pair of vertices $i, j\in \Omega$, the corresponding
components $x_{i}$, $x_{j}$ of $X$ are equal.
\end{lemma}

\begin{definition}{\bf\cite{LSQ}} \label{de02,03} 
Let $r\geq 1$, $G = (V, E)$ be a hypergraph with $u\in V$ and $e_{1}, \ldots, e_{r} \in
E$, such that $u \notin e_{i}$ for $i = 1$, $\ldots$, $r$. Suppose that $v_{i}\in e_{i}$ and write $e^{'}_{i}
= (e_{i}\setminus \{v_{i} \})\cup\{u\}$ $(i = 1, \ldots, r)$. Let $G^{'} = (V, E^{'})$ be the hypergraph with $E^{'} = (E\setminus \{e_{1}$, $\ldots$, $e_{r}\})\cup
\{e^{'}_{1}, \ldots, e^{'}_{r}\}$. Then we say that $G^{'}$ is obtained from $G$ by moving edges $(e_{1}, \ldots, e_{r})$
from $(v_{1}, \ldots, v_{r})$ to $u$.
\end{definition}

\begin{lemma}{\bf \cite{LSQ}} \label{le02,04} 
Let $r\geq 1$, $G$ be a connected uniform hypergraph, $G^{'}$ be the hypergraph obtained
from $G$ by moving edges $(e_{1}$, $\ldots$, $e_{r})$ from $(v_{1}, \ldots, v_{r})$ to $u$, and $G^{'}$ contain no
multiple edges. If in the principal eigenvector $X$ of $G$,
$x_{u}\geq \max_{1\leq i\leq r} \{x_{v_{i}}\}$, then $\rho(G^{'}
) > \rho(G)$.
\end{lemma}

Let $G = (V, E)$ be a $k$-uniform hypergraph, and $e = \{u_{1}$, $u_{2}$, $\ldots$, $u_{k}\}$, $f = \{v_{1}, v_{2}, \ldots, v_{k}\}$
be two edges of $G$, and let $e^{'} = (e \setminus U_{1})\cup V_{1}$, $f^{'} = (f \setminus V_{1}) \cup U_{1}$, where $U_{1} = \{u_{1}$, $u_{2}$, $\ldots$, $u_{r}\}$,
$V_{1} = \{v_{1}$, $v_{2}$, $\ldots$, $v_{r}\}$, $1\leq r\leq k-1$. Let $G^{'} = (V, E^{'})$ be the hypergraph with $E^{'} =
(E \setminus\{e, f\})\cup \{e^{'}, f^{'}\}$. Then we say that $G^{'}$ is obtained from $G$ by $e\rightleftharpoons
{u_{1},u_{2},\ldots,u_{r}\atop
{v_{1},v_{2},\ldots,v_{r}}} f$ or $e\rightleftharpoons
{U_{1}\atop
{V_{1}}}f $.
Let $G = (V, E)$ be a connected $k$-uniform hypergraph, and $X$ be an eigenvector of $G$. For
the simplicity of the notation, we write: $x_{S} =\prod_{v\in S} x_{v}$ where $S\subseteq V$, and for an edge $e$, we write $x_{e} =\prod_{v\in e} x_{v}$.

\begin{lemma}{\bf \cite{PXLW}}\label{le02,05} 
Let $G = (V, E)$ be a connected $k$-uniform hypergraph, and $e = \{u_{1}$, $u_{2}$, $\ldots$, $u_{k}\}$, $f = \{v_{1}$, $v_{2}$, $\ldots$, $v_{k}\}$
be two edges of $G$, $U_{1} = \{u_{1}$, $u_{2}$, $\ldots$, $u_{r}\}$,
$V_{1} = \{v_{1}$, $v_{2}$, $\ldots$, $v_{r}\}$, $1\leq r\leq k-1$. $G^{'}$ is obtained from $G$ by $e\rightleftharpoons
{U_{1}\atop{V_{1}}}f $. Let $X$ be the principal eigenvector of G, $U_{2} = e \setminus U_{1}$, $V_{2} = f \setminus V_{1}$. If $x_{U_{1}} \geq x_{V_{1}}$, $x_{U_{2}} \leq x_{V_{2}}$, and $G^{'}$ is connected, then
$\rho(G) \leq \rho(G^{'})$. Moreover, if one of the two inequalities is strict, then $\rho(G) < \rho(G^{'})$.
\end{lemma}

\section{Main results}

\begin{lemma}\label{le03.01} 
Let $G\in \mathcal {H}$, $\rho(G)=\rho_{max}$, $u$ be a $M_{a}$-vertex in $G$. Then $deg(u) \geq 2$.
\end{lemma}

\begin{proof}
We prove this result by contradiction. Suppose $deg(u) =1$. Because $G$ has at least $2$ edges and $G$ is connected, it follows that $u$ is adjacent to a vertex with degree at least 2, say $v$ for convenience. Assume that $u$ and $v$ are in the same edge $e$, and assume that for $1\leq i\leq deg(v)-1$, edge $e_{i}\neq e$ is incident with $v$. We let $e^{'}_{i}=(e_{i}\setminus \{v\})\cup \{u\}$ for $1\leq i\leq deg(v)-1$, and let $G^{'}=G-\sum^{deg(v)-1}_{i=1}e_{i}+\sum^{deg(v)-1}_{i=1}e^{'}_{i}$. Then $\rho(G)< \rho(G^{'})$ by Lemma \ref{le02,04}. It is a contradiction because $G\cong G^{'}$. Thus we get that $deg(u) \geq 2$.
This completes the proof. \ \ \ \ \ $\Box$
\end{proof}

\begin{lemma}\label{le03.02} 
Let $G\in \mathcal {H}$ and $\rho(G)=\rho_{max}$. Then the unique cycle in $G$ is a $2$-cycle, and in the $2$-cycle, there is a $M_{a}$-vertex.
\end{lemma}

\begin{proof}
Let $\mathcal {M}$ be a maximal matching of $G$, $X$ be
the $principal$ $eigenvector$ of $G$, and let $C=v_{1}e_{1}v_{2}e_{2}\cdots v_{q}e_{q}v_{1}$ be the unique cycle in $G$.

We claim that there is a $M_{a}$-vertex in $C$. Otherwise, suppose that there is no $M_{a}$-vertex in $C$, and assume that $u$ is a $M_{a}$-vertex in $G$. Denote by $P$ the shortest path from $u$ to $C$. Suppose that $P\cap C=\{v_{s}\}$. If $s\in [q]$, then without loss of generality, assume that $s=1$ for convenience. Note that at most one of $e_{1}$ and $e_{q}$ is in $\mathcal {M}$. If $e_{1}\notin\mathcal {M}$, then we let $e^{'}_{1}=(e_{1}\setminus \{v_{1}\})\cup\{u\}$, and let $G^{'}=G-e_{1}+e^{'}_{1}$; if $e_{1}\in\mathcal {M}$, then we let $e^{'}_{q}=(e_{q}\setminus \{v_{1}\})\cup\{u\}$, and let $G^{'}=G-e_{q}+e^{'}_{q}$. Then $\mathcal {M}$ is also a matching of $G^{'}$ and $G^{'}\in \mathcal {H}$. Using Lemma \ref{le02,04} gets that $\rho(G)< \rho(G^{'})$, which contradicts that $\rho(G)=\rho_{max}$. If $s\notin [q]$, then $v_{s}\in e_{j}$ for some $1\leq j\leq q$. Without loss of generality, suppose $j=1$ for convenience. Similar to the case that $s\in [q]$, we get a $G^{'}\in \mathcal {H}$ that $\rho(G)< \rho(G^{'})$, which contradicts that $\rho(G)=\rho_{max}$. This means that our claim holds.

By the above claim, suppose $v_{t}\in V(C)$ is a $M_{a}$-vertex. Next, we prove $q=2$. We prove this by contradiction. Suppose $q\geq 3$.

{\bf Case 1} $t\in [q]$. Without loss of generality, assume $t=1$ for convenience. Note that at most one of $e_{1}$ and $e_{q}$ is in $\mathcal {M}$. Without loss of generality, assume that $e_{q}\notin\mathcal {M}$. In $X$, if $x(v_{2})\geq x(v_{q})$, let $e^{'}_{q}=(e_{q}\setminus \{v_{q}\})\cup\{v_{2}\}$, and let $G^{'}=G-e_{q}+e^{'}_{q}$. Then $\mathcal {M}$ is also a matching of $G^{'}$ and $G^{'}\in \mathcal {H}$. Using Lemma \ref{le02,04} gets that $\rho(G)< \rho(G^{'})$, which contradicts that $\rho(G)=\rho_{max}$. In $X$, if $x(v_{2})< x(v_{q})$ and $e_{1}\notin\mathcal {M}$, we let $e^{'}_{1}=(e_{1}\setminus \{v_{2}\})\cup\{v_{q}\}$, and let $G^{'}=G-e_{1}+e^{'}_{1}$; if $x(v_{2})< x(v_{q})$ and $e_{1}\in\mathcal {M}$, then $e_{2}\notin\mathcal {M}$, and then we let $e^{'}_{2}=(e_{2}\setminus \{v_{2}\})\cup\{v_{1}\}$, and let $G^{'}=G-e_{2}+e^{'}_{2}$. Similarly, we get that $G^{'}\in \mathcal {H}$ and $\rho(G)< \rho(G^{'})$, which contradicts that $\rho(G)=\rho_{max}$.

{\bf Case 2} $t\notin [q]$. Then $v_{t}\in e_{j}$ for some $1\leq j\leq q$. Without loss of generality, suppose $j=1$ for convenience. Similar to Case 1, we can get a unicyclic hypergraph $G^{'}\in \mathcal {H}$ that $\rho(G)< \rho(G^{'})$, which contradicts that $\rho(G)=\rho_{max}$.

By Case 1 and Case 2, it follows that $q=2$. Consequently, the lemma follows from the above narrations.
This completes the proof. \ \ \ \ \ $\Box$
\end{proof}

\begin{lemma}\label{le03.03} 
Let $G \in \mathcal {H}$, $\rho(G)=\rho_{max}$, $C=v_{1}e_{1}v_{2}e_{2}v_{1}$ be the unique cycle in $G$, $u\in V(C)$ be a $M_{a}$-vertex. For each vertex $v$ other than $u$, denote by $P_{u,v}$ a shortest path from $u$ to $v$. We have

$\mathrm{(1)}$ if $V(P_{u,v})\cap V(C)=\{u\}$, then  $d(u, v) \leq 2$; moreover, if $v$ is a NPP-vertex, then $d(u, v) =1$;

$\mathrm{(2)}$ if $\|V(P_{u,v})\cap V(C)\|\geq 2$ and $V(P_{u,v})\cap V(C)$ can be contained in exactly one edge of $C$, then  $d(u, v) \leq 2$; moreover, if $v$ is a NPP-vertex, then $d(u, v) =1$;

$\mathrm{(3)}$ if $\|V(P_{u,v})\cap V(C)\|\geq 2$ and $V(P_{u,v})\cap V(C)$ can not be contained in exactly one edge of $C$, then  $d(u, v) \leq 3$; moreover, if $v$ is a NPP-vertex, then $d(u, v) \leq 2$;

$\mathrm{(4)}$ if $u\in\{v_{1}, v_{2}\}$, then  $d(u, v) \leq 2$; moreover, if $v$ is a NPP-vertex, then $d(u, v) =1$.

\end{lemma}

\begin{proof}
Let $\mathcal {M}$ be a maximal matching of $G$.

(1) Under the condition $V(P_{u,v})\cap V(C)=\{u\}$, we prove two claims.

{\bf Claim 1} If $v$ is a PP-vertex, then $d(u, v) \leq 2$. Otherwise, suppose that $v$ is a PP-vertex and $d(u, v) \geq 3$. Assume that $v\in e_{3}$, and assume that $e_{4}$ is the nonpendant edge adjacent to $e_{3}$ which is in the path $P_{u,v}$ and $e_{3}\cap e_{4}=w$. Clearly, $u\notin e_{4}$.

If $e_{3}\notin \mathcal {M}$, then we let $e^{'}_{3}= (e_{3}\setminus\{w\})\cup\{u\}$, $G^{'}=G-e_{3}+e^{'}_{3}$. Then $\mathcal {M}$ also is a matching of $G^{'}$ and $G^{'}\in \mathcal {H}$. Using Lemma \ref{le02,04} gets that $\rho(G)< \rho(G^{'})$, which contradicts that $\rho(G)=\rho_{max}$.

If $e_{3}\in \mathcal {M}$, then $e_{4}\notin \mathcal {M}$. Suppose $e_{5}$ in $P_{u,v}$ is adjacent to $e_{4}$ that $e_{5}\neq e_{3}$ and $e_{5}\cap e_{4}=\{w_{1}\}$. Let $e^{'}_{4}= (e_{4}\setminus\{w_{1}\})\cup\{u\}$, $G^{'}=G-e_{4}+e^{'}_{4}$. Similar to the case that $e_{3}\notin \mathcal {M}$, we get that $G^{'}\in \mathcal {H}$ and $\rho(G)< \rho(G^{'})$, which contradicts that $\rho(G)=\rho_{max}$.

As a result, Claim 1 holds.

{\bf Claim 2} If $v$ is a NPP-vertex, we have $d(u, v) =1$.
Suppose that there is a NPP-vertex $v$ other than $u$ that $d(u, v) \geq 2$. Assume that in $P_{u,v}$, $v$ is in the edge $e_{3}$, $e_{3}$ is adjacent to $e_{4}$ and $e_{3}\cap e_{4}=\{w\}$, where $w\neq v$. From the assumption that $v$ is a NPP-vertex, we know that $e_{3}$ is not a pendant edge. Thus there is a vertex $t$ of $e_{1}$ other than $w$ is adjacent to at least one edge which is not on $P_{u,v}$ ($t=v$ possible). Here, $d(u,t)=d(u,v)$. Suppose that $P^{'}=te_{s_{1}}v_{s_{1}}e_{s_{2}}v_{s_{2}}\cdots e_{s_{j}}v_{s_{j}}$ is a longest path form $t$ that $V(P_{u,v})\cap V(P^{'})=\{t\}$. Note that $L(P^{'})\geq 1$ and the maximality of the length of $P^{'}$ form $t$. Then $e_{s_{j}}$ must be a pendant edge and $v_{s_{j}}$ is a PP-vertex. Otherwise, it contradicts that $P^{'}$ is longest. Note that $G$ is unicyclic.
Thus $d(u,v_{s_{j}})=d(u, t)+L(P^{'})\geq 3$, which contradicts Claim 1 because $v_{s_{j}}$ is a PP-vertex. Then Claim 2 holds.

 As a result, (1) follows. (2), (3) and (4) are proved similar to (1).
  This completes the proof. \ \ \ \ \ $\Box$
\end{proof}

By Lemma \ref{le03.03}, we have the following two Corollaries.

\begin{corollary}\label{le03.04} 
Let $G \in \mathcal {H}$, $\rho(G)=\rho_{max}$, $C=v_{1}e_{1}v_{2}e_{2}v_{1}$ be the unique cycle in $G$, $u\in V(C)$ be a $M_{a}$-vertex. If there exists a vertex $v$ other than $u$ such that $d(u, v) = 3$, then  $u\notin\{v_{1}, v_{2}\}$, $v$ is a PP-vertex, $\|P_{u,v}\cap C\|\geq 2$ and $P_{u,v}\cap C$ can not be contained in exactly one edge of $C$, where $P_{u,v}$ is a shortest path from $u$ to $v$.
\end{corollary}

\begin{corollary}\label{le03.05} 
Let $G \in \mathcal {H}$, $\rho(G)=\rho_{max}$, $C=v_{1}e_{1}v_{2}e_{2}v_{1}$ be the unique cycle in $G$, $u\in V(C)$ be a $M_{a}$-vertex. For each vertex $v$ other than $u$, denote by $P_{u,v}$ a shortest path from $u$ to $v$. We have

$\mathrm{(1)}$ if $V(P_{u,v})\cap V(C)=\{u\}$, $d(u, v) = 2$, then $v$ is a PP-vertex;

$\mathrm{(2)}$ if $\|V(P_{u,v})\cap V(C)\|\geq 2$, $P_{u,v}\cap C$ can be contained in exactly one edge of $C$, and $d(u, v) =2$, then $v$ is a PP-vertex;

$\mathrm{(3)}$ if $u\in\{v_{1}, v_{2}\}$ and  $d(u, v) =2$, then $v$ is a PP-vertex.
\end{corollary}

\begin{lemma}{\bf \cite{YQY}}\label{le03.06} 
For a tensor $\mathcal {T} \geq 0$, we have
$$c=\min_{1\leq i\leq n}\sum_{i_{2},\ldots,i_{m}=1}
\mathcal {T}_{ii_{2},\ldots,i_{m}}
\leq \rho(\mathcal {T}) \leq \max_{1\leq i\leq n}\sum_{i_{2},\ldots,i_{m}=1}
\mathcal {T}_{ii_{2},\ldots,i_{m}}=C.$$

\end{lemma}

\begin{lemma}{\bf \cite{MFYZ}} \label{le03.07} 
Suppose that tensor $\mathcal {T} \geq 0$ is weakly irreducible. Then either equality in Lemma \ref{le03.06} holds if and only if $c=C$.
\end{lemma}

From Lemmas \ref{le03.06} and \ref{le03.07}, for a uniform connected hypergraph, we have $\delta\leq \rho(G)\leq \Delta$ with equality if and only if $\delta=\Delta$.

\begin{lemma}\label{le03.08} 
Let $G \in \mathcal {H}$, $\rho(G)=\rho_{max}$, $C=v_{1}e_{1}v_{2}e_{2}v_{1}$ be the unique cycle in $G$, $u\in V(C)$ be a $M_{a}$-vertex. Then $u\in\{v_{1}$, $v_{2}\}$, i.e. except $v_{1}$, $v_{2}$, no other vertex of $C$ is $M_{a}$-vertex.
\end{lemma}

\begin{proof}
We prove this lemma by contradiction. Let $X$ be the principal vector of $G$. Suppose $u\in e_{1}$ where $u\neq v_{1}$, $u\neq v_{2}$. Assume that $e_{1}=\{v_{1}$, $v_{2}$, $u$, $v_{1,1}$, $v_{1,2}$, $\cdots$, $v_{1,k-3}\}$, $e_{2}=\{v_{1}$, $v_{2}$, $v_{2,1}$, $v_{2,2}$, $\cdots$, $v_{2,k-2}\}$. By Lemma \ref{le03.01}, we know that $deg(u)\geq 2$. Note that $G$ is unicyclic. Thus, there is an edge $e_{3}$ incident with $u$ that $e_{3}\cap C=\{u\}$, where we denote by $e_{3}=\{u$, $v_{3,1}$, $v_{3,2}$, $\cdots$, $v_{3,k-1}\}$.

{\bf Case 1} There is some $1\leq j\leq k-2$ that $deg(v_{2,j})\geq 2$. Note that $G$ is unicyclic. Assume that $e_{4}=\{v_{2,j}$, $v_{4,1}$, $v_{4,2}$, $\cdots$, $v_{4,k-1}\}$ is an edge incident with $v_{2,j}$ that $e_{4}\cap C=\{v_{2,j}\}$. Note that $d(u,v_{4,i})=3$ for $1\leq i\leq k-1$. By Corollary \ref{le03.04}, it follows that $e_{4}$ is a pendant edge. We claim that there is a  maximal matching $\mathcal {M}$ of $G$ which contains no edge $e_{2}$. Otherwise, suppose $M$ is a maximal matching of $G$ which contains edge $e_{2}$. Then $M$ contains no edge $e_{4}$. Then letting $\mathcal {M}=(M\setminus \{e_{2}\})\cup\{e_{4}\}$ makes our claim holding. Now, let $e^{'}_{2}=(e_{2}\setminus \{v_{2}\})\cup \{u\}$, $G^{'}=G-e_{2}+e^{'}_{2}$. Note that $\mathcal {M}$ is also a matching of $G^{'}$ and $G^{'}\in \mathcal {H}$. Using Lemma \ref{le02,04} gets that $\rho(G)< \rho(G^{'})$, which contradicts $\rho(G)=\rho_{max}$.

{\bf Case 2} For $1\leq i\leq k-2$ that $deg(v_{2,i})= 1$. By Lemma \ref{le02,02}, then for $1\leq i< j\leq k-2$, we have $x(v_{2,i})=x(v_{2,j})$. Let $\mathcal {M}$ be a maximal matching of $G$.

If $x(v_{3,k-1})\geq x(v_{2})$, then let $e^{'}_{1}=(e_{1}\setminus \{v_{2}\})\cup \{v_{3,k-1}\}$, $G^{'}=G-e_{1}+e^{'}_{1}$. Note that at most one of $e_{1}$ and $e_{3}$ is in $\mathcal {M}$. If $e_{1}\notin \mathcal {M}$, then $\mathcal {M}$ is also a maximal matching of $G^{'}$; if $e_{1}\in \mathcal {M}$, then $(\mathcal {M}\setminus \{e_{1}\})\cup \{e^{'}_{1}\}$ is a maximal matching of $G^{'}$. Thus $G^{'}\in \mathcal {H}$. Using Lemma \ref{le02,04} gets that $\rho(G)< \rho(G^{'})$, which contradicts $\rho(G)=\rho_{max}$.

If $x(v_{3,k-1})< x(v_{2})$ and $x(u)\prod^{k-2}_{i=1}x(v_{3,i})\geq x(v_{1})\prod^{k-2}_{i=1}x(v_{2,i})$, then let $G^{'}$ be obtained by $e_{2}\rightleftharpoons
{v_{2}\atop
{v_{3,k-1}}} e_{3}$. Denote by $e^{'}_{2}=(e_{2}\setminus \{v_{2}\})\cup \{v_{3,k-1}\}$, $e^{'}_{3}=(e_{3}\setminus \{v_{3,k-1}\})\cup \{v_{2}\}$. Then we can get a matching $\mathcal {M}^{'}$ of $G^{'}$ from $\mathcal {M}$ by replacing $e_{2}$ with $e^{'}_{2}$ if $e_{2}\in \mathcal {M}$ and replacing $e_{3}$ with $e^{'}_{3}$ if $e_{3}\in \mathcal {M}$. Thus $\alpha(G^{'})\geq \alpha(G)$ and $G^{'}\in \mathcal {H}$. Using Lemma \ref{le02,05} gets that $\rho(G)< \rho(G^{'})$, which contradicts $\rho(G)=\rho_{max}$.

If $x(v_{3,k-1})< x(v_{2})$ and $x(u)\prod^{k-2}_{i=1}x(v_{3,i})< x(v_{1})\prod^{k-2}_{i=1}x(v_{2,i})$, then $1\leq\frac{x(u)}{x(v_{1})}< \frac{\prod^{k-2}_{i=1}x(v_{2,i})}{\prod^{k-2}_{i=1}x(v_{3,i})}$. It follows that $\prod^{k-2}_{i=1}x(v_{3,i})< \prod^{k-2}_{i=1}x(v_{2,i})$. Let $G^{'}$ be obtained by $e_{2}\rightleftharpoons
{v_{2}, v_{2,1}, v_{2,2}, \cdots, v_{2,k-2}\atop
{v_{3,1}, v_{3,2}, \cdots, v_{3,k-1}}} e_{3}$. Similar to the case that $x(v_{3,k-1})< x(v_{2})$ and $x(u)\prod^{k-2}_{i=1}x(v_{3,i})\geq x(v_{1})\prod^{k-2}_{i=1}x(v_{2,i})$, we get that $\alpha(G^{'})\geq \alpha(G)$ and $G^{'}\in \mathcal {H}$. Using Lemma \ref{le02,05} gets that $\rho(G)< \rho(G^{'})$, which contradicts $\rho(G)=\rho_{max}$.

By above narrations, we get that except $v_{1}$, $v_{2}$, no other vertex of $C$ is a $M_{a}$-vertex. Using Lemma \ref{le03.02} gets that $u\in\{v_{1}$, $v_{2}\}$.
This completes the proof. \ \ \ \ \ $\Box$
\end{proof}

Furthermore, by Corollary \ref{le03.05} and Lemma \ref{le03.08}, we have the following Corollary \ref{le03.09}.

\begin{corollary}\label{le03.09} 
Let $G \in \mathcal {H}$, $\rho(G)=\rho_{max}$, $C=v_{1}e_{1}v_{2}e_{2}v_{1}$ be the unique cycle in $G$, $u\in\{v_{1}$, $v_{2}\}$ be a $M_{a}$-vertex. Then every edge not incident with $u$ must be a pendant edge.
\end{corollary}

\begin{lemma}\label{le03.10} 
Let $G \in \mathcal {H}$, $\rho(G)=\rho_{max}$, $C=v_{1}e_{1}v_{2}e_{2}v_{1}$ be the unique cycle, $v_{1}$ be a $M_{a}$-vertex in $G$. Then

$\mathrm{(1)}$ for each vertex $v$ other than $v_{1}$ and $v_{2}$,  we have $deg(v) \leq 2$; moreover, every such vertex $v$ with degree $2$ must be incident with exactly one pendant edge;

$\mathrm{(2)}$ for vertex $v_{2}$,  we have $deg(v_{2}) \leq 3$; moreover, if $deg(v_{2})=3$, then $v_{2}$ must be incident with exactly one pendant edge.
\end{lemma}

\begin{proof}
By Corollary \ref{le03.05} and Lemma \ref{le03.08}, for every vertex $v$ other than $v_{1}$ with $deg(v) \geq 2$, it follows that the distance $d(v_{1},v)\leq 1$. Thus for vertex $v$ with $deg(v) \geq 2$ other than $v_{1}$ and $v_{2}$, we suppose that both $v_{1}$, $v$ are in edge $e_{i}$, and $deg(v)=s \geq 3$. Other than $e_{i}$, assume that $e_{i_{1}}$, $e_{i_{2}}$, $\ldots$, $e_{i_{s-1}}$ are incident with $v$. By Corollary \ref{le03.09}, it follows that $e_{i_{1}}$, $e_{i_{2}}$, $\ldots$, $e_{i_{s-1}}$ are all pendant edges here. Let $\mathcal {M}$ be a maximal matching of $G$. Note that at most one  $e_{i_{j}}$ ($1\leq j\leq s-1$) is in $\mathcal {M}$. Without loss of generality, suppose one of $e_{i}$, $e_{i_{1}}$ is in $\mathcal {M}$. Then $e_{i_{j}}\notin \mathcal {M}$ for $j=2$, $\ldots$, $s-1$. Let $e^{'}_{i_{j}}=(e_{i_{j}}\setminus \{v\})\cup\{v_{1}\}$ for $2\leq j\leq s-1$, $G^{'}=G-\sum_{j=2}^{s-1}e_{i_{j}}+\sum_{j=2}^{s-1}e^{'}_{i_{j}}$. Then $\mathcal {M}$ is also a matching of $G^{'}$. Thus $\alpha(G^{'})\geq\alpha(G)$ and $G^{'}\in \mathcal {H}$. Using Lemma \ref{le02,04} gets that $\rho(G)< \rho(G^{'})$, which contradicts that $\rho(G)=\rho_{max}$. This means that $deg(v) \leq 2$ for each vertex $v$ other than $v_{1}$ and $v_{2}$. Combined with Corollary \ref{le03.09}, $\mathrm{(1)}$ follows.
$\mathrm{(2)}$ is gotten similar to $\mathrm{(1)}$.
This completes the proof. \ \ \ \ \ $\Box$
\end{proof}

Let $G \in \mathcal {H}$, $C=v_{1}e_{1}v_{2}e_{2}v_{1}$ be the unique cycle in $G$, where $e_{1}=\{v_{1}$, $v_{2}$, $v_{1,1}$, $v_{1,2}$, $\cdots$, $v_{1,k-2}\}$, $e_{2}=\{v_{1}$, $v_{2}$, $v_{2,1}$, $v_{2,2}$, $\cdots$, $v_{2,k-2}\}$. $C$ is called a $v_{1}$-$complex$ cycle if $C$ satisfies that no pendant edges is incident with $v_{1}$, and satisfies one of the following (1) and (2):

(1) $deg(v_{2})\geq 3$;

(2) at least one in $\{v_{1,1}$, $v_{1,2}$, $\cdots$, $v_{1,k-2}\}$ is with degree more than $1$, and at least one in $\{v_{2,1}$, $v_{2,2}$, $\cdots$, $v_{2,k-2}\}$ is with degree more than $1$.

\begin{lemma}\label{le03.11} 
Let $G\in \mathcal {H}$ and $\rho(G)=\rho_{max}$, $C=v_{1}e_{1}v_{2}e_{2}v_{1}$ be the unique cycle, $v_{1}$ be a $M_{a}$-vertex in $G$. For an edge $e$ not incident with $v_{1}$, let $e^{'}=(e\setminus N_{G}(v_{1}))\cup\{v_{1}\}$, and $G^{'}=G-e+e^{'}$. Then

$\mathrm{(1)}$ $C$ is not a $v_{1}$-complex cycle;

$\mathrm{(2)}$ $\rho(G)< \rho(G^{'})$;

$\mathrm{(3)}$ $\alpha(G^{'})=\alpha(G)-1$.

\end{lemma}

\begin{proof}
Using Lemma \ref{le02,04}, we get $\rho(G)< \rho(G^{'})$. By Corollary \ref{le03.09}, for edge $e$ not incident with $v_{1}$, we know that $e$ is a pendant edge. Noting that $G$ is unicyclic, by Lemma \ref{le03.03}, it follows that $\|e\cap N_{G}(v_{1})\|=1$. Suppose that $e\cap N_{G}(v_{1})=\{u\}$. Then by Lemma \ref{le03.10}, it follows that $e$ is the unique pendant edge incident with $u$. Let $\mathcal {M}$ be a maximal matching of $G$. Then $\mathcal {M}\setminus\{e\}$ is a matching of $G^{'}$. For any matching $M$ of $G^{'}$, if $e^{'}\notin M$, then $M$ is also a matching of $G$; if $e^{'}\in M$, noting that in $G^{'}$, any edge incident with $v_{1}$ other than $e^{'}$ is not in $M$, then $(M\setminus \{e^{'}\})\cup\{e\}$ is a matching of $G$.
Thus $\alpha(G)-1\leq \alpha(G^{'})\leq \alpha(G)$.

Now, we prove $C$ is not a $v_{1}$-complex cycle. Ohterwise, suppose $C$ is a $v_{1}$-complex cycle. We assert that there is a maximal matching $\mathcal {M}^{'}$ of $G$ which contains no any edge incident with $v_{1}$. For a maximal matching $\mathcal {M}$ of $G$, if $\mathcal {M}$ contains no any edge incident with $v_{1}$, then we let $\mathcal {M}^{'}=\mathcal {M}$.
If there is an edge, say $e_{s}$, which is incident with $v_{1}$, satisfying that $e_{s}\in \mathcal {M}$, then $\mathcal {M}$ contains no any other edge incident with $v_{1}$. Note that $e_{s}$ is not a pendant edge, and note that $C$ is a $v_{1}$-complex cycle. Combining Corollary \ref{le03.09}, we know that $e_{s}$ is adjacent to at least one pendant edge. Then any pendant edge adjacent to $e_{s}$ is not in $\mathcal {M}$. Suppose $e_{i}$ is a pendant edge adjacent to $e_{s}$. Then $\mathcal {M}^{'}=(\mathcal {M}\setminus e_{s})\cup \{e_{i}\}$ satisfies our assertion.
Note that the supposition that $C$ is a $v_{1}$-complex cycle. Then there is an edge not incident with $v_{1}$. For an edge $e$ not incident with $v_{1}$, by Corollary \ref{le03.09}, we know that $e$ is a pendant edge. Here, for edge $e$, if $e\notin \mathcal {M}^{'}$, then $\mathcal {M}^{'}$ is also a matching of $G^{'}$; if $e\in \mathcal {M}^{'}$, then $(\mathcal {M}^{'}\setminus \{e\})\cup\{e^{'}\}$ is a matching of $G^{'}$.
Thus $\alpha(G)\leq \alpha(G^{'})$ and $G^{'}\in \mathcal {H}$, but $\rho(G)< \rho(G^{'})$, which contradicts that $\rho(G)=\rho_{max}$. As a result, it follows that $C$ is not a $v_{1}$-complex cycle.

Next, we prove $\alpha(G^{'})=\alpha(G)-1$ by contradiction. Suppose $ \alpha(G^{'})=\alpha(G)$. Note that $C$ is not a $v_{1}$-complex cycle in $G$, and note that any edge not incident with $v_{1}$ is a pendant edge by Corollary \ref{le03.09}. Thus for proving $\alpha(G^{'})=\alpha(G)-1$, we need consider two following cases.

{\bf Case 1} There are some pendant edges incident with $v_{1}$ in $G$. Assume that $e_{j}$ is a pendant edge incident with $v_{1}$ in $G$. Then $e_{j}$ is also a pendant edge incident with $v_{1}$ in $G^{'}$ and $e_{j}\neq e^{'}$. We assert that there is a maximal matching $\mathcal {M}^{\circ}$ of $G^{'}$ which contains $e_{j}$.

Suppose $\mathbb{M}$ is a maximal matching of $G^{'}$. We claim that $\mathbb{M}$ contains an edge incident with $v_{1}$. Otherwise, in $G^{'}$, if $\mathbb{M}$ contains no any edge incident with $v_{1}$, then $\mathbb{M}\cup\{e^{'}\}$ is a matching with larger cardinality than $\mathbb{M}$, which contradicts the maximality of $\mathbb{M}$. This ensures our claim holding.

If $\mathbb{M}$ does not contain any pendant edge incident with $v_{1}$, then  a nonpendant edge incident with $v_{1}$, say $e_{w}$, is contained in $\mathbb{M}$. We let $\mathcal {M}^{\circ}=(\mathbb{M}\setminus e_{w})\cup\{e_{j}\}$. If $e^{'}\in \mathbb{M}$, then we let $\mathcal {M}^{\circ}=(\mathbb{M}\setminus e^{'})\cup\{e_{j}\}$. Thus $\mathcal {M}^{\circ}$ is also a maximal matching of $G^{'}$, which satisfies our assertion.

Note that $\mathcal {M}^{\circ}$ is a matching. It follows the fact that among all edges incident with $v_{1}$, $e_{j}$ is the exactly one in $\mathcal {M}^{\circ}$. Note that edge $e$ not incident with $v_{1}$ is a pendant edge. By Lemma \ref{le03.03}, it follows that in  $G$, pendant edge $e$ is adjacent to a nonpendant edge $e_{r}$, where $v_{1}\in e_{r}$. Note that $e_{r}\notin \mathcal {M}^{\circ}$. Then $\mathcal {M}^{\circ}\cup \{e\}$ is a matching of $G$ with cardinality $\alpha(G^{'})+1$. Note the supposition that $\alpha(G)= \alpha(G^{'})$. Thus $\mathcal {M}^{\circ}\cup \{e\}$ contradicts the matching number of $G$. As a result, it follows that $\alpha(G^{'})=\alpha(G)-1$.

{\bf Case 2} In $G$, no pendant edge is incident with $v_{1}$. Note that $C$ is not a $v_{1}$-complex cycle. Then $deg(v_{2})= 2$, and no pendant edge is adjacent to one of $e_{1}$, $e_{2}$. Suppose that no pendant edge is adjacent to $e_{1}$ in $G$. As Case 1, we can prove that there is a maximal matching $\mathcal {M}^{\circ}$ of $G^{'}$ which contains $e_{1}$, and prove that $\mathcal {M}^{\circ}\cup \{e\}$ is a matching of $G$ with cardinality $\alpha(G)+1$, which contradicts the matching number of $G$. Thus $\alpha(G^{'})=\alpha(G)-1$.

To sum up, the lemma follows as desired.
This completes the proof. \ \ \ \ \ $\Box$
\end{proof}

From Lemma \ref{le03.11}, we get the following Corollary \ref{le03.12} immediately.

\begin{corollary}\label{le03.12} 
Let $G\in \mathcal {H}$ and $\rho(G)=\rho_{max}$, $C=v_{1}e_{1}v_{2}e_{2}v_{1}$ be the unique cycle, $v_{1}$ be a $M_{a}$-vertex in $G$. If $deg(v_{2})=3$, then there must be a pendant edge incident with $v_{1}$.
\end{corollary}

\begin{lemma}\label{le03.13} 
Let $G\in \mathcal {H}$ and $\rho(G)=\rho_{max}$, $C=v_{1}e_{1}v_{2}e_{2}v_{1}$ be the unique cycle, $v_{1}$ be a $M_{a}$-vertex in $G$. Then any maximal matching of $G$ consists of exactly one pendant edge incident with $v_{1}$, or exactly one edge of $C$ in which each vertex other than $v_{1}$ is not incident with any pendant edge, and all pendant edges not incident with $v_{1}$.
\end{lemma}

\begin{proof}
Let $\mathcal {M}$ be a maximal matching of $G$. By Lemma \ref{le03.11}, we know that $C$ is not a $v_{1}$-complex cycle. For proving this lemma, we need consider two cases next.

{\bf Case 1} There are some pendant edges incident with $v_{1}$ in $G$. As the proof for the claim about $\mathbb{M}$ in Case 1 of Lemma \ref{le03.11}, we get that there must be an edge incident with $v_{1}$ is in $\mathcal {M}$. Suppose that an edge $e_{3}$ incident with $v_{1}$ is in $\mathcal {M}$. Then any other edge incident with $v_{1}$ is not in $\mathcal {M}$.

We claim that $e_{3}$ is a pendant edge or an edge of $C$ in which each vertex other than $v_{1}$ is not incident with any pendant edge. We prove this claim by contradiction. Suppose this claim does not hold. Then two subcases need consider.

{\bf Subcase 1.1} $e_{3}$ is a nonpendant edge which is not in $C$.

By Corollary \ref{le03.09}, there must be a pendant edge $e_{4}$ adjacent to $e_{3}$ where $e_{3}\cap e_{4}=\{v_{3}\}$, $v_{3}\neq v_{1}$. Because $e_{3}\in \mathcal {M}$, then $e_{4} \notin \mathcal {M}$. Let $e^{'}_{4}=(e_{4}\setminus \{v_{3}\})\cup \{v_{1}\}$, $G^{'}=G-e_{4}+e^{'}_{4}$. Then $\mathcal {M}$ is also a matching of $G^{'}$ and $G^{'}\in \mathcal {H}$. Using Lemma \ref{le02,04} gets $\rho(G^{'})>\rho(G)$, which contradicts $\rho(G)=\rho_{max}$.

{\bf Subcase 1.2} $e_{3}$ is an edge in $C$ which is adjacent to some pendant edge $e_{t}$, where $e_{t}$ is not incident with $v_{1}$. For this case, as Subcase 1.1, we can get a hypergraph $G^{'}\in \mathcal {H}$ that $\rho(G^{'})>\rho(G)$, which contradicts $\rho(G)=\rho_{max}$.

By Subcase 1.1 and Subcase 1.2, our claim holds.

For a pendant edge $e_{5}$ not incident with $v_{1}$, by Lemma \ref{le03.03} and Corollary \ref{le03.05}, $e_{5}$ is adjacent to an edge $e_{6}$ where $v_{1}\in e_{6}$, $e_{6}\neq e_{3}$. Note that $G$ is unicyclic. Therefore, $\|e_{5}\cap e_{6}\|=1$. Suppose $e_{5}\cap e_{6}=\{v_{5}\}$ where $v_{5}\neq v_{1}$. By Lemma \ref{le03.10}, $e_{5}$ is the unique pendant edge incident with $v_{5}$. If $e_{5}\notin \mathcal {M}$, noting that $e_{3}\in \mathcal {M}$, $e_{6}\notin \mathcal {M}$, we get that $\{e_{5}\}\cup \mathcal {M}$ is also a matching of $G$, which contradicts the maximality of $\mathcal {M}$. Thus $e_{5}\in \mathcal {M}$. Then it follows that the lemma holds for Case 1.

{\bf Case 2} In $G$, no pendant edge is incident with $v_{1}$, $deg(v_{2})= 2$, and no pendant edge is adjacent to one of $e_{1}$, $e_{2}$. For this case, as Case 1, it is proved that the lemma holds.

This completes the proof. \ \ \ \ \ $\Box$
\end{proof}

Furthermore, by Lemma \ref{le03.13}, we get the following Corollary \ref{le03.14}.

\begin{corollary}\label{le03.14} 
Let $G\in \mathcal {H}$ and $\rho(G)=\rho_{max}$, $C=v_{1}e_{1}v_{2}e_{2}v_{1}$ be the unique cycle, $v_{1}$ be a $M_{a}$-vertex in $G$. If there is a pendant edge incident with $v_{1}$, then
there is a maximal matching of $G$ consists of exactly one pendant edge incident with $v_{1}$, and all other pendant edges not incident with $v_{1}$.
\end{corollary}

Let $\mathbb{F}$ be a connected $k$-uniform hypergraph ($k\geq 3$) with at least $2$ edges and $u\in V(\mathbb{F})$, $e=\{v_{1}$, $v_{2}$, $\ldots$, $v_{k-1}$, $u\}$ be a nonpendant edge incident with $u$, and $e_{1}$, $e_{2}$, $\ldots$, $e_{t}$ ($t\leq k-1$) be the pendant edges incident with vertices $v_{1}$, $v_{2}$, $\ldots$, $v_{t}$ of $e$ respectively, where $t\leq k-1$, $e_{i}\cap e=\{v_{i}\}$, $deg(v_{i})=2$ for $1\leq i\leq t$, $deg(v_{i})=1$ for $t+1\leq i\leq k-1$ if $t\leq k-2$. Here, letting $t=0$ means that $deg(v_{i})=1$ for $1\leq i\leq k-1$.

\begin{lemma}\label{le03.15} 
For hypergraph $\mathbb{F}$, let $\rho=\rho(\mathbb{F})$. Then in the principal eigenvector $X$ of $\mathbb{F}$, $x(v_{i})=\frac{x_{u}}{\rho(1-\frac{1}{{\rho}^{k}})^{\frac{t+1}{k}}}$ for $1\leq i\leq t$, and $x(v_{i})=\frac{x(u)}{\rho(1-\frac{1}{{\rho}^{k}})^{\frac{t}{k}}}$ for $t+1\leq i\leq k-1$ if $t\leq k-2$. Moreover, if $t=k-1$, then $x(v_{i})=\frac{x(u)}{\rho-\frac{1}{\rho^{k-1}}}$ for $1\leq i\leq t$; if $t=0$, then $x(v_{i})=\frac{x(u)}{\rho}$ for $1\leq i\leq k-1$.
\end{lemma}

\begin{proof}
Note that $\mathbb{F}$ is connected and $\mathbb{F}$ has at least 2 edges. By Lemma \ref{le03.06} and Lemma \ref{le03.07}, it follows that $\rho> 1$. Suppose $e_{i}=\{v_{i}$, $v_{i_{1}}$, $v_{i_{2}}$, $\ldots$, $v_{i_{k-1}}\}$ for $1\leq i\leq t$. In the principal eigenvector $X$ of $\mathbb{F}$, let $x(v_{1})=w$, and let $x(v_{t+1})=f$ if $t\leq k-2$. By Lemma \ref{le02,02}, then $x(v_{i})=w$, $x(v_{i_{j}})=\frac{w}{\rho}$ for $1\leq i\leq t$, $1\leq j\leq k-1$; $x(v_{i})=f$ for $t+1\leq i\leq k-1$. Let $s=k-1-t$. If $s\geq 1$, then
$$\left \{\begin{array}{ll}
               \rho f^{k-1}=f^{s-1}w^{t}x(u),  \ & \ \ \ \ \ (1)\\
              \rho w^{k-1}=f^{s}w^{t-1}x(u)+(\frac{w}{\rho})^{k-1}  \ & \ \ \ \ \ (2).
             \end{array}\right.$$
It follows that $w=\frac{x_{u}}{\rho(1-\frac{1}{{\rho}^{k}})^{\frac{t+1}{k}}}$, $f=(1-\frac{1}{\rho^{k}})^{\frac{1}{k}}w=\frac{x_{u}}{\rho(1-\frac{1}{{\rho}^{k}})^{\frac{t}{k}}}$.

If $s=0$, then $x(v_{i})=w=\frac{x_{u}}{\rho-\frac{1}{\rho^{k-1}}}$ for $1\leq i\leq t$ follows from (2); if $t=0$, then $x(v_{i})=f=\frac{x(u)}{\rho}$ for $1\leq i\leq k-1$ follow from (1).
Thus the lemma follows as desired.
This completes the proof. \ \ \ \ \ $\Box$
\end{proof}

\setlength{\unitlength}{0.5pt}
\begin{center}
\begin{picture}(671,256)
\put(522,95){\circle*{4}}
\put(506,111){\circle*{4}}
\put(102,188){\circle*{4}}
\put(36,161){\circle*{4}}
\put(204,237){\circle*{4}}
\put(156,83){\circle*{4}}
\put(175,223){\circle*{4}}
\put(228,175){\circle*{4}}
\put(142,206){\circle*{4}}
\put(196,157){\circle*{4}}
\put(102,110){\circle*{4}}
\put(46,67){\circle*{4}}
\put(132,95){\circle*{4}}
\put(74,47){\circle*{4}}
\put(84,124){\circle*{4}}
\put(72,173){\circle*{4}}
\put(66,137){\circle*{4}}
\put(266,206){\circle*{4}}
\put(225,74){\circle*{4}}
\put(152,140){\circle*{4}}
\put(605,231){\circle*{4}}
\put(560,67){\circle*{4}}
\put(571,209){\circle*{4}}
\put(547,196){\circle*{4}}
\put(602,150){\circle*{4}}
\put(461,153){\circle*{4}}
\put(492,124){\circle*{4}}
\put(438,203){\circle*{4}}
\put(490,163){\circle*{4}}
\put(512,172){\circle*{4}}
\put(568,126){\circle*{4}}
\put(671,206){\circle*{4}}
\put(638,173){\circle*{4}}
\put(432,86){\circle*{4}}
\put(7,95){\circle*{4}}
\put(541,79){\circle*{4}}
\put(497,40){\circle*{4}}
\qbezier(36,161)(104,233)(204,237)
\qbezier(36,161)(154,165)(204,237)
\qbezier(36,161)(107,139)(156,83)
\qbezier(36,161)(69,94)(156,83)
\qbezier(175,223)(229,204)(228,175)
\qbezier(175,223)(202,175)(228,175)
\qbezier(142,206)(189,186)(196,157)
\qbezier(142,206)(166,160)(196,157)
\qbezier(102,110)(49,96)(46,67)
\qbezier(102,110)(64,57)(46,67)
\qbezier(132,95)(76,72)(74,47)
\qbezier(132,95)(99,41)(74,47)
\qbezier(204,237)(249,238)(266,206)
\qbezier(204,237)(238,204)(266,206)
\qbezier(156,83)(204,91)(225,74)
\qbezier(156,83)(197,57)(225,74)
\qbezier(102,188)(147,173)(152,140)
\qbezier(102,188)(127,135)(152,140)
\qbezier(547,196)(596,185)(602,150)
\qbezier(547,196)(570,151)(602,150)
\qbezier(438,203)(466,204)(490,163)
\qbezier(438,203)(444,178)(490,163)
\qbezier(512,172)(565,156)(568,126)
\qbezier(605,231)(665,232)(671,206)
\qbezier(605,231)(654,196)(671,206)
\qbezier(571,209)(631,201)(638,173)
\qbezier(571,209)(610,171)(638,173)
\qbezier(512,172)(543,119)(568,126)
\qbezier(461,153)(476,91)(560,67)
\qbezier(461,153)(559,158)(605,231)
\qbezier(461,153)(534,227)(605,231)
\qbezier(461,153)(532,127)(560,67)
\qbezier(461,153)(426,114)(432,86)
\qbezier(461,153)(459,95)(432,86)
\qbezier(36,161)(29,106)(7,95)
\qbezier(36,161)(0,123)(7,95)
\qbezier(541,79)(514,35)(497,40)
\qbezier(497,40)(496,64)(541,79)
\put(322,134){$\Longrightarrow$}
\put(123,15){$\mathbb{K}$}
\put(530,15){$\mathbb{K}^{'}$}
\put(250,-15){Fig. 3.1. $\mathbb{K}$ and $\mathbb{K}^{'}$}
\put(24,164){$u$}
\put(449,159){$u$}
\put(194,245){$v_{1,1}$}
\put(156,230){$v_{1,2}$}
\put(86,195){$v_{1,4}$}
\put(43,176){$v_{1,5}$}
\put(591,239){$v_{1,1}$}
\put(548,216){$v_{1,2}$}
\put(140,71){$v_{2,1}$}
\put(103,111){$v_{2,3}$}
\put(548,55){$v_{2,1}$}
\put(541,80){$v_{2,2}$}
\put(132,96){$v_{2,2}$}
\put(500,179){$v_{1,4}$}
\put(489,152){$v_{1,5}$}
\put(469,63){\circle*{4}}
\qbezier(522,95)(473,83)(469,63)
\qbezier(469,63)(490,55)(522,95)
\put(520,100){$v_{2,3}$}
\end{picture}
\end{center}

Let $\mathbb{K}$ be a connected $k$-uniform hypergraph ($k\geq 3$) with at least $2$ edges and $u\in V(\mathbb{K})$, $e_{1}=\{v_{1,1}$, $v_{1,2}$, $\ldots$, $v_{1,k-1}$, $u\}$ be a nonpendant edge incident with $u$, and $e_{1,1}$, $e_{1,2}$, $\ldots$, $e_{1,s}$ ($1\leq s\leq k-2$) be the pendant edges incident with vertices $v_{1,1}$, $v_{1,2}$, $\ldots$, $v_{1,s}$ of $e_{1}$ respectively, where $e_{1,i}\cap e_{1}=\{v_{1,i}\}$, $deg(v_{1,i})=2$ for $1\leq i\leq s$, $deg(v_{1,i})=1$ for $s+1\leq i\leq k-1$; $e_{2}=\{v_{2,1}$, $v_{2,2}$, $\ldots$, $v_{2,k-1}$, $u\}$ be aother nonpendant edge incident with $u$, and $e_{2,1}$, $e_{2,2}$, $\ldots$, $e_{2,t}$ ($1\leq t\leq s$) be the pendant edges incident with vertices $v_{2,1}$, $v_{2,2}$, $\ldots$, $v_{2,t}$ of $e_{2}$ respectively, where $e_{2,i}\cap e_{2}=\{v_{2,i}\}$, $deg(v_{2,i})=2$ for $1\leq i\leq t$, $deg(v_{2,i})=1$ for $t+1\leq i\leq k-1$. Let $e^{'}_{2,t}=(e_{2,t}\setminus \{v_{2,t}\})\cup \{v_{1,s+1}\}$, $\mathbb{K}^{'} =\mathbb{K}-e_{2,t}+e^{'}_{2,t}$ (for example, see Fig. 3.1).

\begin{lemma}\label{le03.16} 
For hypergraph $\mathbb{K}$, we have $\rho(\mathbb{K}^{'})> \rho(\mathbb{K})$ and $\alpha(\mathbb{K}^{'})\geq\alpha(\mathbb{K})$.
\end{lemma}

\begin{proof}
Note that $\mathbb{K}$ is connected and $\mathbb{K}$ has at least 2 edges. By Lemma \ref{le03.06} and Lemma \ref{le03.07}, it follows that $\rho(\mathbb{K})> 1$. Let $X$ be the principal eigenvector of $\mathbb{K}$. Combining Lemma \ref{le02,02}, we let $x(v_{1,i})=w_{1}$ for $1\leq i\leq s$; $x(v_{1,i})=f_{1}$ for $s+1\leq i\leq k-1$; $x(v_{2,i})=w_{2}$ for $1\leq i\leq t$; $x(v_{2,i})=f_{2}$ for $t+1\leq i\leq k-1$. Then by Lemma \ref{le03.15}, $w_{1}\geq w_{2}$, $f_{1}\geq f_{2}$, $w_{1}> f_{1}$, $w_{2}> f_{2}$. Thus $x(u)\prod^{s}_{i=1}x(v_{1,i})\prod^{k-1}_{i=s+2}x(v_{1,i})> x(u)\prod^{t-1}_{i=1}x(v_{2,i})\prod^{k-1}_{i=t+1}x(v_{2,i})$. If $f_{1}\geq w_{2}$, then $\rho(\mathbb{K}^{'})> \rho(\mathbb{K})$ by Lemma \ref{le02,04}. If $f_{1}< w_{2}$, then $\mathbb{K}^{'}$ is obtained from $\mathbb{K}$ by $e_{1}\rightleftharpoons
{v_{1,s+1}\atop{v_{2,t}}}e_{2}$, and then $\rho(\mathbb{K}^{'})> \rho(\mathbb{K})$ by Lemma \ref{le02,05}.

As Lemma \ref{le03.13}, we get that hypergraph $\mathbb{K}$ has a maximal matching $\mathcal {M}$ which contains no $e_{1}$ and $e_{2}$, but contains all $e_{1,i}$ for $1\leq i\leq s$ and contains all $e_{2,i}$ for $1\leq i\leq t$. Then $(\mathcal {M}\setminus\{e_{2,t}\})\cup \{e^{'}_{2,t}\}$ is a matching of $\mathbb{K}^{'}$. Thus $\alpha(\mathbb{K}^{'})\geq \alpha(\mathbb{K})$.
This completes the proof. \ \ \ \ \ $\Box$
\end{proof}

\begin{lemma}\label{le03.17} 
Let $G\in \mathcal {H}$ and $\rho(G)=\rho_{max}$, $C=v_{1}e_{1}v_{2}e_{2}v_{1}$ be the unique cycle, $v_{1}$ be a $M_{a}$-vertex in $G$. If no pendant edge is incident with any vertex in $e_{1}\setminus \{v_{1}\}$, or no pendant edge is incident with any vertex in $e_{2}\setminus \{v_{1}\}$, then except $e_{1}$ and $e_{2}$, any other edge incident with $v_{1}$ is a pendant edge.
\end{lemma}

\begin{proof}
Suppose no pendant edge is incident with any vertex in $e_{1}\setminus \{v_{1}\}$ in $G$. Denote by $e_{1}=\{v_{1}$, $v_{1,1}$, $v_{1,2}$, $\ldots$, $v_{1,k-2}$, $v_{2}\}$. We prove this lemma by contradiction. Suppose this lemma does not hold, and assume that $e_{3}$ is a nonpendant edge incident with $v_{1}$ where $e_{3}\neq e_{1}$, $e_{3}\neq e_{2}$. Denote by $e_{3}=\{v_{1}$, $v_{3,1}$, $v_{3,2}$, $\ldots$, $v_{3,k-2}$, $v_{3,k-1}\}$. By Corollary \ref{le03.09} and Lemma \ref{le03.10}, assume that $e_{3,k-1}$, $e_{3,k-2}$, $\ldots$, $e_{3,k-s}$ ($1\leq s\leq k-1$) are the pendant edges incident with vertices $v_{3,k-1}$, $v_{3,k-2}$, $\ldots$, $v_{3,k-s}$ of $e_{3}$ respectively. Note that $G$ is unicyclic. Then by Lemmas \ref{le03.06} and \ref{le03.07}, it follows that $\rho(G)> 1$. Let $X$ be the principal eigenvector of $G$. By Lemma \ref{le03.15}, then $x(v_{1,1})=x(v_{1,2})=\cdots=x(v_{1,k-2})$, $x(v_{3,k-1})=x(v_{3,k-2})=\cdots=x(v_{3,k-s})$, $x(v_{3,k-s+1})=x(v_{3,k-s+2})=\cdots=x(v_{3,1})$ and $x(v_{3,1})< x(v_{3,k-1})$ if $s\leq k-2$. Let $\mathcal {M}$ be a maximal matching of $G$. Then by Lemma \ref{le03.13}, we get that $e_{3,k-i}\in \mathcal {M}$ for $1\leq i\leq s$, but $e_{3}\notin \mathcal {M}$; if there is a pendant edge $e$ incident with a vertex in $e_{2}\setminus \{v_{1}\}$, then $e\in \mathcal {M}$ and $e_{2}\notin \mathcal {M}$. Note that at most one of $e_{1}$ and $e_{2}$ is in $\mathcal {M}$. Thus without loss of generality, we suppose that $e_{2}\notin \mathcal {M}$.

{\bf Case 1} $x(v_{3,k-1})\geq x(v_{2})$. Then let $e^{'}_{2}=(e_{2}\setminus\{v_{2}\})\cup\{v_{3,k-1}\}$, and let $G^{'} =G -e_{2}+e^{'}_{2}$. Then $\mathcal {M}$ is also a maximal matching of $G^{'}$ and $G^{'}\in \mathcal {H}$. Using Lemma \ref{le02,04} gets $\rho(G^{'})> \rho(G)$, which contradicts $\rho(G)=\rho_{max}$.

{\bf Case 2} $x(v_{3,k-1})< x(v_{2})$.

{\bf Subcase 2.1} $x(v_{3,k-1})\leq x(v_{1,1})$. Denote by $e^{'}_{3,k-i}=(e_{3,k-i}\setminus \{v_{3,k-i}\})\cup\{v_{1,k-i}\}$ for $2\leq i\leq s$, and $e^{'}_{3,k-1}=(e_{3,k-1}\setminus \{v_{3,k-1}\})\cup\{v_{2}\}$. Let
$G^{'} =G -e_{3,k-1}+e^{'}_{3,k-1}-\sum^{s}_{i=2}e_{3,k-i}+\sum^{s}_{i=2}e^{'}_{3,k-i}$, and let $U_{i}=e_{3,k-i}\setminus \{v_{3,k-i}\}$ for $1\leq i\leq s$. Note that if $e_{1}\notin \mathcal {M}$, then $\mathcal {M}$ is also a maximal matching of $G^{'}$; if $e_{1}\in \mathcal {M}$, then $\mathcal {M}^{'}=(\mathcal {M}\setminus \{e_{1}\})\cup\{e_{3}\}$ is a matching of $G^{'}$. Thus $G^{'}\in \mathcal {H}$.
Note that $$\rho(G^{'})-\rho(G)\geq X^{T}\mathcal {A}(G^{'})X^{k-1}-X^{T}\mathcal {A}(G)X^{k-1}\hspace{7.1cm}$$$$=k(x_{e^{'}_{3,k-1}}-x_{e_{3,k-1}})+k(\sum^{s}_{i=2}x_{e^{'}_{3,k-i}}-\sum^{s}_{i=2}x_{e_{3,k-i}})\hspace{2.5cm}$$
$$=kx_{U_{1}}(x(v_{2})-x(v_{3,k-1}))+k\sum^{s}_{i=2}x_{U_{i}}(x(v_{1,k-i})-x(v_{3,k-i}))>0.\hspace{0.3cm}$$ Therefore, $\rho(G^{'})> \rho(G)$ by Lemma \ref{le01.02}, which contradicts $\rho(G)=\rho_{max}$.

{\bf Subcase 2.2} $x(v_{3,k-1})> x(v_{1,1})$.

{\bf Subcase 2.2.1} $x(v_{1,1})x(v_{1,2})\cdots x(v_{1,k-2})\geq x(v_{3,1})x(v_{3,2})\cdots x(v_{3,k-2})$. Then $$x^{k-2}(v_{1,1})\geq x^{k-1-s}(v_{3,1})x^{s-1}(v_{3,k-s})\Longrightarrow \frac{x^{k-1-s}(v_{1,1})}{x^{k-1-s}(v_{3,1})}\geq \frac{x^{s-1}(v_{3,k-s})}{x^{s-1}(v_{1,1})}=\frac{x^{s-1}(v_{3,k-1})}{x^{s-1}(v_{1,1})}>1.$$

Let $Q_{1}=\{v_{1,k-s}$, $v_{1,k-s+1}$, $\ldots$, $v_{1,k-2}\}$, $Q_{2}=\{v_{3,k-s}$, $v_{3,k-s+1}$, $\ldots$, $v_{3,k-2}\}$, $Q_{3}=e_{1}\setminus Q_{1}$, $Q_{4}=e_{3}\setminus Q_{2}$, $e^{'}_{1}=Q_{3}\cup Q_{2}$, $e^{'}_{3}=Q_{4}\cup Q_{1}$, $U=e_{3,k-1}\setminus \{v_{3,k-1}\}$ and $e^{'}_{3,k-1}=U\cup\{v_{2}\}$. Then $x_{Q_{3}}>x_{Q_{4}}$ and $x_{Q_{2}}>x_{Q_{1}}$. Let
$G^{'} =G -e_{3,k-1}-e_{1}-e_{2}+e^{'}_{3,k-1}+e^{'}_{1}+e^{'}_{2}$. Note that if $e_{1}\notin \mathcal {M}$, then $\mathcal {M}$ is also a maximal matching of $G^{'}$; if $e_{1}\in \mathcal {M}$, then $\mathcal {M}^{'}=(\mathcal {M}\setminus \{e_{1}\})\cup\{e_{3}\}$ is a matching of $G^{'}$. Thus $G^{'}\in \mathcal {H}$. Note that  $$\rho(G^{'})-\rho(G)\geq X^{T}\mathcal {A}(G^{'})X^{k-1}-X^{T}\mathcal {A}(G)X^{k-1}=k(x_{Q_{3}}-x_{Q_{4}})(x_{Q_{2}}-x_{Q_{1}})+kx_{U}(x(v_{2})-x(v_{3,k-1}))>0.$$ Then $\rho(G^{'})> \rho(G)$ by Lemma \ref{le01.02}, which contradicts $\rho(G)=\rho_{max}$.

{\bf Subcase 2.2.2} $x(v_{1,1})x(v_{1,2})\cdots x(v_{1,k-2})< x(v_{3,1})x(v_{3,2})\cdots x(v_{3,k-2})$. Let $Q_{1}=\{v_{1,1}$, $v_{1,2}$, $\ldots$, $v_{1,k-2}\}$, $Q_{2}=\{v_{3,1}$, $v_{3,2}$, $\ldots$, $v_{3,k-2}\}$, $e^{'}_{1}=(e_{1}\setminus Q_{1})\cup Q_{2}$, $e^{'}_{3}=(e_{3}\setminus Q_{2})\cup Q_{1}$, and $e^{'}_{3,k-1}=(e_{3,k-1}\setminus \{v_{3,k-1}\})\cup\{v_{2}\}$. Let
$G^{'} =G -e_{3,k-1}-e_{1}-e_{2}+e^{'}_{3,k-1}+e^{'}_{1}+e^{'}_{2}$. As Subcase 2.2.1, we get that $\alpha(G^{'})\geq \alpha(G)$, $G^{'}\in \mathcal {H}$ and $\rho(G^{'})> \rho(G)$, which contradicts $\rho(G)=\rho_{max}$.

By the above narrations, we get that except $e_{1}$ and $e_{2}$,  no other nonpendant edge is incident with $v_{1}$. Then the lemma follows. This completes the proof. \ \ \ \ \ $\Box$
\end{proof}

Let $\mathbb{W}$ be a connected $k$-uniform hypergraph ($k\geq 3$) with at least $2$ edges, and $C=v_{1}e_{1}v_{2}e_{2}v_{1}$ be a cycle in $\mathbb{W}$. Denote by $e_{1}=\{v_{1}$, $v_{1,1}$, $v_{1,2}$, $\ldots$, $v_{1,k-2}$, $v_{2}\}$, and let $e_{1,k-2}$, $e_{1,k-3}$, $\ldots$, $e_{1,k-t}$ ($2\leq t\leq k-1$) be the pendant edges incident with vertices $v_{1,k-2}$, $v_{1,k-3}$, $\ldots$, $v_{1,k-t}$ of $e_{1}$ respectively, where $e_{1,k-i}\cap e_{1}=\{v_{1,k-i}\}$, $deg(v_{1,i})=2$ for $2\leq i\leq t$, $deg(v_{1,i})=1$ for $t+1\leq i\leq k-1$ if $t\leq k-2$. Here, letting $t=1$ means that $deg(v_{1,k-i})=1$ for $2\leq i\leq k-1$. Similar to Lemma \ref{le03.15}, for hypergraph $\mathbb{W}$, we have the following Lemma \ref{le03.18}.

\begin{lemma}\label{le03.18} 
For graph $\mathbb{W}$, let $\rho=\rho(\mathbb{W})$. Then in the principal eigenvector $X$ of $\mathbb{W}$, $x(v_{1,k-i})=(\frac{x(v_{1})x(v_{2})}{\varphi})^{\frac{1}{2}}$ for $2\leq i\leq t$, and $x(v_{i})=(1-\frac{1}{{\rho}^{k}})^{\frac{1}{k}}(\frac{x(v_{1})x(v_{2})}{\varphi})^{\frac{1}{2}}$ for $t+1\leq i\leq k-1$ if $t\leq k-2$, where $\varphi=\rho(1-\frac{1}{{\rho}^{k}})^{\frac{t+1}{k}}$. Moreover, if $t=k-1$, then $x(v_{1,k-i})=\left (\frac{x(v_{1})x(v_{2})}{\rho-\frac{1}{\rho^{k-1}}}\right)^{\frac{1}{2}}$ for $2\leq i\leq t$; if $t=1$, then $x(v_{1,k-i})=(\frac{x(v_{1})x(v_{2})}{\rho})^{\frac{1}{2}}$ for $2\leq i\leq k-1$.
\end{lemma}

\begin{lemma}\label{le03.19} 
Let $G\in \mathcal {H}$ and $\rho(G)=\rho_{max}$, $C=v_{1}e_{1}v_{2}e_{2}v_{1}$ be the unique cycle in $G$, $v_{1}$ be a $M_{a}$-vertex, and $X$ be the principal eigenvector of $G$.

$\mathrm{(1)}$ If there is a nonpendant edge incident with $v_{1}$ which is neither of $e_{1}$ and $e_{2}$, then $deg(v_{2})=3$, and $x(v_{2})>x(v)$ for any $v\in (V(G)\setminus\{v_{1}, v_{2}\})$.

$\mathrm{(2)}$ Assume that there is a pendant edge $e_{l}$ incident with a vertex in $e_{1}\setminus\{v_{1}\}$, and there is a pendant edge $e_{p}$ incident with a vertex in $e_{2}\setminus\{v_{1}\}$ ($e_{l}=e_{p}$ possible, and $e_{l}=e_{p}$ means that $e_{l}$ is incident with $v_{2}$). Then $deg(v_{2})=3$, and $x(v_{2})>x(v)$ for any $v\in (V(G)\setminus\{v_{1}, v_{2}\})$.

$\mathrm{(3)}$ Assume that there is a pendant edge $e_{l}$ incident with a vertex in $e_{1}\setminus\{v_{1}\}$, and there is a pendant edge incident with $v_{1}$. Then $deg(v_{2})=3$, and $x(v_{2})>x(v)$ for any $v\in (V(G)\setminus\{v_{1}, v_{2}\})$.
\end{lemma}

\begin{proof}
We prove (1) first. Note the condition that there is a nonpendant edge incident with $v_{1}$ which is neither of $e_{1}$ and $e_{2}$. By Lemma \ref{le03.17}, it follows that both $e_{1}$ and $e_{2}$ are adjacent to some pendant edges respectively, where these pendant edges are not incident with $v_{1}$. Let $\mathcal {M}$ be a maximal matching of $G$. Denote by $e_{1}=\{v_{1}$, $v_{1,1}$, $v_{1,2}$, $\ldots$, $v_{1,k-2}$, $v_{2}\}$, $e_{2}=\{v_{1}$, $v_{2,1}$, $v_{2,2}$, $\ldots$, $v_{2,k-2}$, $v_{2}\}$. Combining Corollary \ref{le03.09} and Lemma \ref{le03.10}, suppose $e_{1,k-2}$, $e_{1,k-3}$, $\ldots$, $e_{1,k-s}$ ($2\leq s\leq k-1$) are the pendant edges incident with vertices $v_{1,k-2}$, $v_{1,k-3}$, $\ldots$, $v_{1,k-s}$ of $e_{1}$ respectively; $e_{2,k-2}$, $e_{2,k-3}$, $\ldots$, $e_{2,k-t}$ ($2\leq t\leq k-1$) are the pendant edges incident with vertices $v_{2,k-2}$, $v_{2,k-3}$, $\ldots$, $v_{2,k-t}$ of $e_{2}$ respectively. Note that $G$ is unicyclic. Then by Lemmas \ref{le03.06} and \ref{le03.07}, it follows that $\rho(G)> 1$. By Lemma \ref{le02,02} and Lemma \ref{le03.18}, we get that
$x(v_{1,k-2})=x(v_{1,k-3})=\ldots=x(v_{1,k-s})$, $x(v_{1,k-s-1})=x(v_{1,k-s-2})=\ldots=x(v_{1,1})$ and $x(v_{1,k-2})>x(v_{1,1})$ if $s\leq k-2$; $x(v_{2,k-2})=x(v_{2,k-3})=\ldots=x(v_{2,k-t})$, $x(v_{2,k-t-1})=x(v_{2,k-t-2})=\ldots=x(v_{2,1})$ and $x(v_{2,k-2})>x(v_{2,1})$ if $t\leq k-2$. By Lemma \ref{le03.15}, it follows that $x(v_{1,k-i})>x(v_{i_{j}})$ where $v_{i_{j}}\in (e_{1,k-i}\setminus \{v_{1,k-i}\})$ and $2\leq i\leq s$; $x(v_{2,k-i})>x(v_{i_{j}})$ where $v_{i_{j}}\in (e_{2,k-i}\setminus \{v_{2,k-i}\})$ and $2\leq i\leq t$. Let $A=V(C)\cup(\cup^{s}_{i=2} e_{1,k-i})\cup(\cup^{t}_{i=2} e_{2,k-i})$, $B=A\setminus\{v_{1}, v_{2}\}$.
Without loss of generality, suppose that $x(v_{1,k-2})\geq x(v_{2,k-2})$. If $x(v_{2})\leq x(v_{1,k-2})$, then let $e^{'}_{2}=(e_{2}\setminus\{v_{2}\})\cup \{v_{1,k-2}\}$, $G^{'} =G -e_{2}+e^{'}_{2}$. Note that $e_{1}\notin \mathcal {M}$ and $e_{2}\notin \mathcal {M}$ by Lemma \ref{le03.13}. Then $\mathcal {M}$ is also a maximal matching of $G^{'}$ and $G^{'}\in \mathcal {H}$. Using Lemma \ref{le02,04} gets that $\rho(G^{'})> \rho(G)$, which contradicts $\rho(G)=\rho_{max}$. This means that $x(v_{2})> x(v_{1,k-2})$. As a result, it follows that $x(v_{2})>x(v)$ for any $v\in B$. Note that if $deg(v_{2})=2$, then we let $e^{'}_{1,k-2}=(e_{1,k-2}\setminus \{v_{1,k-2}\})\cup\{v_{2}\}$, and $G^{''}=G-e_{1,k-2}+e^{'}_{1,k-2}$. Similarly, we get that  $G^{''}\in \mathcal {H}$ and $\rho(G^{''})> \rho(G)$, which contradicts $\rho(G)=\rho_{max}$. This implies that $deg(v_{2})>2$. Combining this and Lemma \ref{le03.10} gets that $deg(v_{2})=3$. Suppose $e_{\mu}$ is the pendant edge incident with $v_{2}$. By Lemma \ref{le03.15}, it follows that $x(v_{2})>x(v_{\mu_{j}})$ where $v_{\mu_{j}}\in (e_{\mu}\setminus \{v_{2}\})$. Let $D=B\cup (e_{\mu}\setminus \{v_{2}\})$. Consequently, we get that $x(v_{2})>x(v)$ for any $v\in D$.

Suppose $e_{f}=\{v_{1}$, $v_{f,1}$, $v_{f,2}$, $\ldots$, $v_{f,k-1}\}$ is an edge incident with $v_{1}$ which is neither of $e_{1}$ and $e_{2}$. If $e_{f}$ is a nonpendant edge, combining Corollary \ref{le03.09} and Lemma \ref{le03.10}, then we assume that $e_{f,k-1}$, $e_{f,k-2}$, $\ldots$, $e_{f,k-\varphi}$ ($1\leq \varphi\leq k-1$) are the pendant edges incident with vertices $v_{f,k-1}$, $v_{f,k-2}$, $\ldots$, $v_{f,k-\varphi}$ of $e_{f}$ respectively. By Lemma \ref{le02,02} and Lemma \ref{le03.15}, if $e_{f}$ is a nonpendant edge, then $x(v_{f,k-1})=x(v_{f,k-2})=\ldots=x(v_{f,k-\varphi})$, $x(v_{f,k-\varphi-1})=x(v_{f,k-\varphi-2})=\ldots=x(v_{f,1})$ and $x(v_{f,k-1})>x(v_{f,1})$ if $\varphi\leq k-2$, and  $x(v_{f,k-i})>x(v_{i_{j}})$ where $v_{i_{j}}\in (e_{f,k-i}\setminus \{v_{f,k-i}\})$ and $1\leq i\leq \varphi$; if $e_{f}$ is a pendant edge, then $x(v_{f,k-1})=x(v_{f,k-2})=\ldots=x(v_{f,1})$. As proved for $x(v_{2})> x(v_{1,k-2})$ in the first paragraph, it follows that $x(v_{f,k-1})<x(v_{2})$. Let $H=e_{f}\cup(\cup^{\varphi}_{i=1} e_{f,k-i})$, $L=H\setminus\{v_{1}\}$ if $e_{f}$ is a nonpendant edge; $L=e_{f}\setminus\{v_{1}\}$ if $e_{f}$ is a pendant edge. Thus it follows that $x(v_{2})>x(v)$ for any $v\in L$.

Notice the arbitrariness of $e_{f}$ and combine that $x(v_{2})>x(v)$ for any $v\in D$.
Then (1) follows as desired.
Combining Lemma \ref{le03.10} and the proof of (1) gets (2); combining Corollary \ref{le03.14} and the proof of (1) gets (3).
This completes the proof. \ \ \ \ \ $\Box$
\end{proof}

\begin{lemma}\label{le03.20} 
Let $G\in \mathcal {H}$ and $\rho(G)=\rho_{max}$, $C=v_{1}e_{1}v_{2}e_{2}v_{1}$ be the unique cycle, $v_{1}$ be a $M_{a}$-vertex in $G$. If there is a nonpendant edge incident with $v_{1}$ which is neither of $e_{1}$ and $e_{2}$, then for each vertex $v\in (V(C)\setminus\{v_{1}\})$, it is incident with exactly one pendant edge.
\end{lemma}

\begin{proof}
Note the condition that there is a nonpendant edge incident with $v_{1}$ which is neither of $e_{1}$ and $e_{2}$. By Lemma \ref{le03.17}, it follows that both $e_{1}$ and $e_{2}$ are adjacent to some pendant edges respectively, where these pendant edges are not incident with $v_{1}$. By Lemma \ref{le03.19}, it follows that there is exactly one pendant edge incident with $v_{2}$. Let $X$ be the principal eigenvector of $G$, and $\mathcal {M}$ be a maximal matching of $G$. Denote by $e_{1}=\{v_{1}$, $v_{1,1}$, $v_{1,2}$, $\ldots$, $v_{1,k-2}$, $v_{2}\}$. Combining Corollary \ref{le03.09} and Lemma \ref{le03.10}, suppose that $e_{1,k-2}$, $e_{1,k-3}$, $\ldots$, $e_{1,k-t}$ ($2\leq t\leq k-1$) are the pendant edges incident with vertices $v_{1,k-2}$, $v_{1,k-3}$, $\ldots$, $v_{1,k-t}$ of $e_{1}$ respectively.
Suppose $e_{f}=\{v_{1}$, $v_{f,1}$, $v_{f,2}$, $\ldots$, $v_{f,k-1}\}$ is a nonpendant edge incident with $v_{1}$ which is neither of $e_{1}$ and $e_{2}$, and combining Corollary \ref{le03.09} and Lemma \ref{le03.10} again, suppose that $e_{f,k-1}$, $e_{f,k-2}$, $\ldots$, $e_{f,k-s}$ ($1\leq s\leq k-1$) are the pendant edges incident with vertices $v_{f,k-1}$, $v_{f,k-2}$, $\ldots$, $v_{f,k-s}$ of $e_{f}$ respectively. By Lemma \ref{le03.19}, we get that $x(v_{2})>x(v)$ for any $v\in (V(G)\setminus\{v_{1}, v_{2}\})$.

Next, we prove that $t= k-1$. Otherwise, suppose that $t\leq k-2$. Then we consider two cases as follows.

{\bf Case 1} $x(v_{1,1})\geq x(v_{f,k-1})$. We let $e^{'}_{f,k-1}=(e_{f,k-1}\setminus \{v_{f,k-1}\})\cup \{v_{1,1}\}$, $G^{'} =G -e_{f,k-1}+e^{'}_{f,k-1}$. Note that $e_{1}\notin \mathcal {M}$ and $e_{f}\notin \mathcal {M}$ by Lemma \ref{le03.13}. Then $\mathcal {M}$ is also a maximal matching of $G^{'}$ and $G^{'}\in \mathcal {H}$. Using Lemma \ref{le02,04} gets that $\rho(G^{'})> \rho(G)$, which contradicts $\rho(G)=\rho_{max}$.

{\bf Case 2} $x(v_{1,1})< x(v_{f,k-1})$.

{\bf Subcase 2.1} $x(v_{1,k-2})x(v_{1,k-3})\cdots x(v_{1,1})\geq x(v_{f,k-1})x(v_{f,k-2})\cdots x(v_{f,2})$. Then $$\frac{x(v_{1,k-2})x(v_{1,k-3})\cdots x(v_{1,2})}{x(v_{f,k-2})x(v_{f,k-3})\cdots x(v_{f,2})}\geq \frac{x(v_{f,k-1})}{x(v_{1,1})}> 1.$$ Note that $x(v_{2})>x(v_{f,1})$. Then $$x(v_{2})x(v_{1,k-2})x(v_{1,k-3})\cdots x(v_{1,2})x(v_{1})>x(v_{f,k-2})x(v_{f,k-3})\cdots x(v_{f,2})x(v_{f,1})x(v_{1}).$$ Let
$G^{'}$ be obtained by $e_{1}\rightleftharpoons
{v_{1,1}\atop
{v_{f,k-1}}}e_{f} $.
 Note that $e_{1}\notin \mathcal {M}$ and $e_{f}\notin \mathcal {M}$ by Lemma \ref{le03.13}. Then $\mathcal {M}$ is also a maximal matching of $G^{'}$ and $G^{'}\in \mathcal {H}$. Using Lemma \ref{le02,05} gets that $\rho(G^{'})> \rho(G)$, which contradicts $\rho(G)=\rho_{max}$.

{\bf Subcase 2.2} $x(v_{1,k-2})x(v_{1,k-3})\cdots x(v_{1,1})< x(v_{f,k-1})x(v_{f,k-2})\cdots x(v_{f,2})$. Let $U_{1}=\{v_{1,k-2}$, $v_{1,k-3}$, $\cdots$, $v_{1,1}\}$, $U_{2}=\{v_{f,k-1}$, $v_{f,k-2}$, $\cdots$, $v_{f,2}\}$, $G^{'}$ be obtained by $e_{1}\rightleftharpoons
{U_{1}\atop
{U_{2}}}e_{f} $. Note that $x(v_{2})>x(v_{f,1})$. As Subcase 2.1, we get that $\mathcal {M}$ is also a maximal matching of $G^{'}$, $G^{'}\in \mathcal {H}$ and $\rho(G^{'})> \rho(G)$, which contradicts $\rho(G)=\rho_{max}$.

By the above narrations, it follows that $t= k-1$. Similarly, for $e_{2}$, we get the same conclusion. Thus the lemma follows as desired.
This completes the proof. \ \ \ \ \ $\Box$
\end{proof}

\begin{lemma}\label{le03.21} 
Let $G\in \mathcal {H}$ and $\rho(G)=\rho_{max}$, $C=v_{1}e_{1}v_{2}e_{2}v_{1}$ be the unique cycle, $v_{1}$ be a $M_{a}$-vertex in $G$. We have

$\mathrm{(1)}$ $m= z+1$ if and only if $deg(v_{1})=2$, $deg(v_{2})=2$, and no pendant edge is adjacent to at least one of $e_{1}$, $e_{2}$;

$\mathrm{(2)}$ if $m\geq z+2$, then there must be a pendant edge incident with $v_{1}$; moreover, if $z\geq 2$, then $deg(v_{2})=3$.
\end{lemma}

\begin{proof}
(1) {\bf Claim} If $deg(v_{1})\geq 3$, then there must be a pendant edge incident with $v_{1}$. Otherwise, assume that no pendant edge is incident with $v_{1}$. Then there is a nonpendant edge incident with $v_{1}$ which is neither of $e_{1}$ and $e_{2}$. Then by Lemma \ref{le03.20}, it follows that $deg(v_{2})=3$. By Corollary \ref{le03.12}, it follows that there must be a pendant edge incident with $v_{1}$. This is a contradiction. Thus the claim holds.

By Corollary \ref{le03.14}, there is a maximal matching of $G$ which contains one pendant edge incident with $v_{1}$ and all other pendant edges not incident with $v_{1}$. Note that now neither of $e_{1}$, $e_{2}$ is in this matching. Then $m\geq z+2$, which contradicts $m= z+1$. Thus $deg(v_{1})=2$. By Corollary \ref{le03.12}, it follows that $deg(v_{2})=2$. By Lemma \ref{le03.19}, it follows that no pendant edge is adjacent to at least one of $e_{1}$, $e_{2}$.

From the above proof, we get that if $deg(v_{1})=2$, then $deg(v_{2})=2$, and no pendant edge is adjacent to at least one one of $e_{1}$, $e_{2}$. Suppose no pendant edge is adjacent to $e_{1}$. By Lemma \ref{le03.13}, there is a maximal matching $\mathcal {M}$ which contains $e_{1}$ and all pendant edges. Then $m= z+1$. Thus (1) follows. Combining (1), Corollary \ref{le03.12} and Lemma \ref{le03.19} gets that if $m\geq z+2$, then $deg(v_{1})\geq3$. Using the above Claim for proving (1) and using Lemma \ref{le03.19} again get (2).
This completes the proof. \ \ \ \ \ $\Box$
\end{proof}

Let $\mathbb{D}$ be a connected $k$-uniform hypergraph ($k\geq 3$) with at least $2$ edges, and $C=v_{1}e_{1}v_{2}e_{2}v_{1}$ be a cycle in $\mathbb{D}$. Denote by $e_{1}=\{v_{1}$, $v_{1,1}$, $v_{1,2}$, $\ldots$, $v_{1,k-2}$, $v_{2}\}$, and let $e_{1,k-2}$, $e_{1,k-3}$, $\ldots$, $e_{1,k-s}$ ($2\leq s\leq k-2$) be the pendant edges incident with vertices $v_{1,k-2}$, $v_{1,k-3}$, $\ldots$, $v_{1,k-s}$ of $e_{1}$ respectively, where $e_{1,k-i}\cap e_{1}=\{v_{1,k-i}\}$, $deg(v_{1,i})=2$ for $2\leq i\leq s$, $deg(v_{1,i})=1$ for $s+1\leq i\leq k-1$. Denote by $e_{2}=\{v_{1}$, $v_{2,1}$, $v_{2,2}$, $\ldots$, $v_{2,k-2}$, $v_{2}\}$, and let $e_{2,k-2}$, $e_{2,k-3}$, $\ldots$, $e_{2,k-t}$ ($2\leq t\leq s$) be the pendant edges incident with vertices $v_{2,k-2}$, $v_{2,k-3}$, $\ldots$, $v_{2,k-t}$ of $e_{2}$ respectively, where $e_{2,k-i}\cap e_{2}=\{v_{2,k-i}\}$, $deg(v_{2,i})=2$ for $2\leq i\leq t$, $deg(v_{2,i})=1$ for $t+1\leq i\leq k-1$. Let $e^{'}_{2,k-t}=(e_{2,k-t}\setminus \{v_{2,k-t}\})\cup \{v_{1,k-s-1}\}$, $\mathbb{D}^{'} =\mathbb{D}-e_{2,k-t}+e^{'}_{2,k-t}$. Similar to Lemma \ref{le03.16}, for hypergraph $\mathbb{D}$, combining Lemma \ref{le03.18}, we have the following Lemma \ref{le03.22}.

\begin{lemma}\label{le03.22} 
For hypergraph $\mathbb{D}$, we have $\rho(\mathbb{D}^{'})> \rho(\mathbb{D})$ and $\alpha(\mathbb{D}^{'})\geq\alpha(\mathbb{D})$.
\end{lemma}

{\bf Proof of Theorem \ref{th01.03}.}
Note that $G$ is connected. Thus we have $m\geq z+1$.

If $m= z+1$, then (1) follows from \ref{le03.21}. Next, we consider the case that $m\geq z+2$.

By Lemmas \ref{le03.01}, \ref{le03.03}, \ref{le03.08}, \ref{le03.10}, \ref{le03.13}, \ref{le03.17}, \ref{le03.19}-\ref{le03.21}, Corollaries \ref{le03.04}, \ref{le03.05}, \ref{le03.09}, \ref{le03.12}, \ref{le03.14}, and using Lemmas \ref{le03.11}, \ref{le03.16} and \ref{le03.22} repeatedly, it follows that $G$ is isomorphic to a $\mathscr{U}(n,k;f;r,s;t,w)$ where $z=\alpha(\mathscr{U}(n,k;f;r,s;t,w))=f+r+s+t(k-1)+w+1$, $v_{1}$ of $\mathscr{U}(n,k;f;r,s;t,w)$ is a $M_{a}$-vertex, and $\mathscr{U}(n,k;f;r,s;t,w)$ satisfies (2)-(5) respectively.
This completes the proof. \ \ \ \ \ $\Box$

Theorem \ref{th01.04} follows from Theorem \ref{th01.03} as a corollary directly.

\small {

}

\end{document}